\documentclass[a4paper]{amsart}
\usepackage{amssymb,amscd,verbatim,latexsym}
\usepackage{enumerate}
\usepackage{tikz-cd}
\usepackage[all,pdf]{xy}

\usepackage[T1]{fontenc}
\usepackage[american]{babel}
\usepackage[utf8]{inputenc}
\usepackage{csquotes}

\usepackage{lmodern}
\usepackage{microtype}

\usepackage[ocgcolorlinks]{hyperref}
\hypersetup{
  colorlinks = true,
  allcolors = [rgb]{0.15, 0.37, 0.41},
  linktocpage = true}

\usepackage[mathscr]{eucal}

\usepackage{color}

\newcommand{\supp}{\operatorname{supp}}

\newcommand{\Ind}{\operatorname{Ind}}

\newcommand{\tto}{\rightrightarrows}

\newcommand{\DR}{\operatorname{DR}}

\newcommand{\ch}{\operatorname{ch}}

\newcommand{\CC}{\mathbb C}

\newcommand{\RR}{\mathbb R}

\newcommand{\ZZ}{\mathbb Z}

\newcommand{\ind}{\operatorname{ind}}

\newcommand{\Ad}{\operatorname{Ad}}

\newcommand{\rk}{\operatorname{rk}}

\newcommand{\maA}{\mathcal A}
\newcommand{\maB}{\mathcal B}
\newcommand{\maC}{\mathcal C}
\newcommand{\maD}{\mathcal D}
\newcommand{\maE}{\mathcal E}
\newcommand{\maF}{\mathcal F}
\newcommand{\maG}{\mathcal G}
\newcommand{\maH}{\mathcal H}
\newcommand{\maI}{\mathcal I}

\newcommand{\maK}{\mathcal K}

\newcommand{\maM}{\mathcal M}
\newcommand{\maN}{\mathcal N}

\newcommand{\maS}{\mathcal S}
\newcommand{\maT}{\mathcal T}

\newcommand{\maZ}{\mathcal Z}

\newcommand{\bs}{\mathbf{s}}
\newcommand{\bt}{\mathbf{t}}

\newtheorem{theorem}{Theorem}[section]
\newtheorem{thm}[theorem]{Theorem}
\newtheorem{prop}[theorem]{Proposition}
\newtheorem{cor}[theorem]{Corollary}

\newtheorem{lem}[theorem]{Lemma}
\newtheorem{nota}[theorem]{Notations}
\theoremstyle{definition}

\newtheorem{dfn}[theorem]{Definition}
\theoremstyle{remark}
\newtheorem{rem}[theorem]{Remark}

\newtheorem{example}[theorem]{Example}
\newtheorem{ex}[theorem]{Example}

\author[Yu Qiao]{Yu Qiao}
\address{School of Mathematics and Statistics, Shaanxi Normal University,
Xi'an, 710119, Shaanxi, China}
\email{yqiao@snnu.edu.cn}

\author[Bing Kwan So]{Bing Kwan So}
\address{School of Mathematics, Jilin University,
Changchun, 130023, Jilin, China}
\email{bkso@graduate.hku.hk}

\author[Hang Wang]{Hang Wang}
\address{Research Center for Operator Algebras, School of Mathematical Sciences, East China Normal University,
3663 North Zhongshan Rd, Shanghai 200062, China}
\email{wanghang@math.ecnu.edu.cn}

\title[An Invariant higher index pairing]{An invariant higher index pairing under Morita equivalence for proper \'etale groupoids}

\begin{document}

\maketitle

\begin{abstract}
For suitable proper \'etale groupoids, we study a part of the higher index that is invariant under Morita equivalence.
We capture this part by a pairing between cyclic homology classes associated with isotropy groups
and cyclic cohomology classes of the groupoid convolution algebra.
We show that the pushforward is compatible with the maps induced by essential equivalences.
This gives a numerical invariant depending only on the Morita equivalence class of the groupoid.
In the case of a finite group action,
we relate the resulting pairing to the equivariant index and to the Atiyah-Bott fixed-point formula.

\end{abstract}

\begin{small}
\hspace*{0.72cm}
\hspace*{2.8cm}
periodic cyclic (co)homology; index pairing. \\
\hspace*{1.1cm}
\end{small}

\tableofcontents

\section{Introduction}

\subsection{Singular objects as groupoids}
In this paper we study index theory of singular objects presented by smooth \'etale groupoids.

One of the main objects we are interested in is the class of orbifolds.
Here we recall the classical definition of effective orbifolds \cite{Ruan;Book}:
An (effective) orbifold $M$ is a topological space
endowed with an orbifold atlas $\{ (\tilde U _\alpha , G_\alpha , \tilde \phi_\alpha) \}$.

Given an effective orbifold atlas $\{ (\tilde U _\alpha , G_\alpha , \tilde \phi_\alpha) \}$ on $M$,
Pronk and Sibih \cite{PronkThesis, PSTarXiv} constructed explicitly a proper \'etale groupoid $\maG$,
with objects $\maG ^{(0)} := \sqcup _\alpha \tilde U _\alpha$,
and arrows $\maG ^{(1)}$ given by equivalence classes (germs) of coordinate transformations
$\tilde U _\alpha \xleftarrow{\phi _\alpha} \tilde U _\gamma \xrightarrow{\phi _\beta} \tilde U _\beta $.
If one uses a different orbifold atlas, then Moerdijk and Pronk in \cite{MP97}
showed that one obtains a {\em Morita equivalent} groupoid.
As for orbifolds that are not effective, see \cite{PSTarXiv}. 

Conversely, given a proper \'etale groupoid $\maG = \maG ^{(1)} \rightrightarrows \maG ^{(0)}$,
Moerdijk and Pronk considered a $\maG$-sheaf, which endows the set of orbits $\maG ^{(0)} / \maG$ with an orbifold structure \cite{MP97}.
Moreover two {\em Morita equivalent} proper \'etale groupoids define isomorphic orbifolds.

More generally, given a regular foliation $F$ on a manifold, one can construct the holonomy groupoid $\maH$.
One can also fix a local transversal and construct an \'etale groupoid. These are all Morita equivalent \cite{MoerdijkFolBook}.

Recall that the category of Lie groupoids, together with bi-bundles and isomorphism of bi-bundles,
forms a weak 2-category, and a bi-bundle is weakly invertible in the 2-categary
if and only if it is induced by a Morita equivalence.
Lastly, given a groupoid $\maG$, Lerman constructed a 2-functor $B$,
which embeds the weak category of bi-bundles into the strict 2-category of stacks,
by assigning each $\maG$ to the category whose object is the $\maG$ principal bundles \cite{Lerman10}.
We shall not go into further details because we only need Morita equivalence.

To conclude, invariance under Morita equivalence is fundamental
to any invariant constructed out of groupoids presenting a singular space.

\subsection{Motivation: noncommutative geometry and index theory}
Ultimately we attempt to construct index type invariant on some `nice' stacks.

The idea of generalizing topological invariants of manifolds to stacks is not new.
For example, Behrend \cite{Behrend;CohStack} and Behrend-Ginot-Nooh-Xu \cite{Xu;StackBook}
constructed the de Rham and singular cohomology of stacks.
Tu, Xu and Laurent-Gengoux studied the twisted $K$-theory of stacks \cite{TXLG04}
(see also Adam and Ruan \cite{Ruan;OrbifoldTwistedK} for the orbifold case).

Even more has been done for the special case of orbifolds.
Farsi and Seaton in \cite{FSPJM10,FSTAMS10} defined the $\Gamma$-sector
for each orbifold groupoid and arbitrary discrete group $\Gamma$, and studied its de Rham cohomology,
generalizing the stringy cohomology of Chen and Ruan \cite{Ruan;OrbifoldCoh, Ruan;OrbifoldCoh2} (where $\Gamma = \ZZ$),
and also Adam and Ruan \cite{Ruan;OrbifoldTwistedK} (where $\Gamma = \ZZ _2 $).
Then Farsi and Seaton \cite{FarsiSeaton11} defined the $\Gamma$-Euler-Satake number
as the sum of the Euler numbers of the $\Gamma$-sectors.

In another series of papers \cite{Farsi92PAMS, Farsi92}, Farsi wrote an orbifold as a quotient $M = \tilde M / H$,
where $H$ is a compact Lie group acting on $\tilde M$, with finite isotropy groups.
Then she paired the compact-group higher index, traditionally called the equivariant index,
of a transversally elliptic operator
to the cyclic homology class represented by the constant 1 function,
and showed that the pairing computes the Kawasaki index, which is indeed an invariant of an orbifold.

In this paper, motivated by the constructions above,
we turn to construction based on noncommutative geometry of \'etale groupoids presenting our orbifold or stack.
It is well known that whenever two groupoids $\maG, \maG'$ are Morita equivalent,
their (compactly supported) convolution algebras are also Morita equivalent \cite{RenBook80}.
Hence one has the isomorphism of $K$-theory $K _\bullet (\maC^\infty_c (\maG)) \cong K _\bullet (\maC^\infty_c (\maG')) $,
and also that of $K$-homology $K ^\bullet (\maC^\infty_c (\maG)) \cong K ^\bullet (\maC^\infty_c (\maG')) $.
Hence we have a Morita invariant notion of abstract pseudo-differential operators,
namely, as elements in $K ^\bullet (\maC^\infty_c (\maG)) $,
and it is natural to consider the index of such operators.

\subsection{Higher and equivariant indices}
The index theory of a group acting on a manifold (that is presented by the action groupoid) has been widely studied.
We briefly recall it here, as a special case, and to serve as further motivation of our work.

Let $G$ be a discrete group acting properly on a manifold $M$.
Then the quotient $M / G$ is an orbifold. The orbifold groupoid of $M / G$ is the action groupoid $M \rtimes G$.
For a $G$-equivariant elliptic pseudo-differential operator $D$, the index class in the $K$-theory
of the group $C^*$-algebra,
$$\Ind_G(D) \in K _0(C^*(G)),$$
is usually called the higher index.

When $G$ is finite, or more generally compact, one has
$K _0(C^*(G)) \cong R(G)$, the representation ring of $G$.
In the most basic case, suppose $G$ is finite and $M$ is compact.
Any equivariant elliptic pseudo-differential operator $D$ is Fredholm,
and the $G$-action preserves the kernel and co-kernel of $D$.
Hence its higher index is the formal difference
$$\Ind_G(D)=[\ker D]-[\operatorname{coker}D] \in R(G) \cong K _0(C^*(G)).$$
This representation-valued index is traditionally called the equivariant index.
Thus, in this paper, ``higher index'' refers to the general $K$-theory-valued class,
whereas ``equivariant index'' refers to its representation-theoretic form for finite or compact groups.
For a general discrete group, the higher index is defined using a parametrix
or by the deformation construction as in \cite{CRWW23}.

There is notable difficulty when one tries to study Morita invariance of the higher index.
For example, let $G$ be any discrete group and $M$ be any manifold.
Let $G$ act on $M \times G$ by right translation on the $G$ component.
This action is free and indeed the action groupoid $(M \times G) \rtimes G$ is Morita equivalent to $M$.
Hence the index of an elliptic operator $D$ on $M$ lies in $K _0 (\CC) \cong \ZZ$,
while the higher index of a $G$-equivariant elliptic operator $\tilde D$ on $M \times G$
lies in $K _0(C^*(G))$
-- which can be as complicated as one likes.
So even within action groupoids it is not obvious how to compare the index as elements in different groups.
We shall give a more elaborate example in Section 3.1.

One can simplify this problem by taking the Chern-Connes character of $D$ to obtain an element in periodic cyclic cohomology,
$\ch (D) \in H P ^\bullet (\maC ^\infty _c (M \rtimes G))$, and then pair with some classes in periodic cyclic homology.
Higher index pairings first appeared in Atiyah's distribution index \cite{AtiyahBook74}.
According to Kasparov \cite{Kasparov16},
the class $\Ind_G(D)$ is the pushforward of the $G$-equivariant $K$-homology class
of $D$ along the constant map $p: M \to pt$.
Pairing with an arbitrary $f \in \mathbb{C}G ^{\Ad {}}$
that can be identified with $H P _0 (C^* (G))$ (since $G$ is finite),
yields Atiyah's distribution index \cite{AtiyahBook74}
$$\langle f , \ch (p _* [D]) \rangle = [\Ind _G(D)](f) := \operatorname{Tr}(f|_{\ker D})
- \operatorname{Tr}(f|_{\operatorname{coker}D}) \in \mathbb{C}.$$
Carrillo Rouse, Wang, and Wang \cite{CRWW23} constructed the dual of the pushforward map
in the pairing $\langle f , \ch (p _* [D]) \rangle$.
Since then higher index pairings have been studied in many geometric settings.
Heat kernel methods provide the classical local framework \cite{BGVBook04};
Connes and Moscovici related higher indices on covering spaces to cyclic cocycles
coming from group cohomology \cite{CM90}; and Roe developed analogous pairings
in coarse geometry \cite{Roe93}.
For groupoids, Pflaum, Posthuma, and Tang obtained a transverse index formula
localized over the inertia groupoid \cite{PPT15}, while Carrillo Rouse, Wang,
and Wang constructed delocalized pairings on the geometric side of the Baum--Connes assembly map \cite{CRWW23}.
Interestingly, by the Atiyah--Bott formula or the delocalized pairing formula in \cite{CRWW23},
only elements in $G$ with fixed points can possibly have non-trivial pairings,
hence one can eliminate the ambiguity in the free action example above just by considering index pairings.


\subsection{The main research problem}
After taking index pairing we end up with a complex number, which can be compared.
However, it is not obvious how to generalize the full equivariant index to all \'etale groupoids
$\maG = \maG ^{(1)} \rightrightarrows \maG ^{(0)}$.
Instead we use cyclic homology classes arising from the isotropy group $\maG_y^y$
(which is just the isotropy subgroup in the action groupoid case).
Pairing these classes with the Chern--Connes character isolates contributions supported by isotropy,
or equivalently by the relevant part of the inertia groupoid.

Our key observation here is that the isotropy subgroup at a point $y \in \maG ^{(0)}$ is an invariant under Morita equivalence,
in the sense that any essential equivalence $\varphi : \maK \to \maG$ naturally restricts to a group isomorphism
from $\maK _x ^x$ to $\maG _y ^y$, where $y = \varphi ^{(0)} (x)$ (see Proposition \ref{IsotropySubgroup}).
It follows that there is a well defined problem of whether there exists a pairing map
$$ H P _\bullet (\maC^\infty _c (\maG _y ^y)) \times H P ^{\bullet + \dim \maG ^{(0)}} (\maC _c ^\infty (\maG)) \to \CC $$
that is compatible with the induced Morita isomorphisms.
Establishing this compatibility produces a numerical higher
index pairing that depends only on the Morita class of the groupoid presentation.

\subsection{Main results of this paper}
Given a proper \'etale groupoids
$\maG$,
our first objective is to produce some periodic homology classes in $H P _\bullet (\maC^\infty _c (\maG))$.
We are guided by the pairing in Ponge and Wang \cite{PW18} for the action groupoid case,
which involves the periodic cyclic homology $H P _\bullet (\maC^\infty _c (G))$.
For these reasons, as a starting point we consider the periodic cyclic homology $H P _\bullet (\maC^\infty _c (\maG ^y_y))$.
In the action groupoid case, observe that the embedding $f _y : \{ y \} \to M$ is a $(M \rtimes G)^y_y$-equivariant map.
Hence Rouse-Wang-Wang \cite{CRWW23} provides us a pushforward map
$$(f_y)_! : H P _\bullet (\maC^\infty _c (\maG ^y_y)) \to H P _{\bullet + \dim \maG ^{(0)}} (\maC^\infty _c (M \rtimes \maG ^y_y)) .$$
Then we have the extension by 0 map $\varepsilon$,
which maps $H P _\bullet (\maC^\infty _c (M \rtimes \maG ^y_y))$ to $H P _\bullet (\maC^\infty _c (M \rtimes G))$.
In section 4, we generalize the construction above to general proper \'etale groupoids that satisfy additional mild assumptions
(recall that locally a proper \'etale groupoid is an action groupoid).

Secondly we must prove our main Theorem \ref{CD6}, which states
that the map $\varepsilon \circ (f_y)_!$ is compatible with the isomorphism induces by essential equivalence.
There we need a more explicit way to describe the isomorphism
$H P _\bullet (\maC^\infty _c (\maK)) \cong H P _\bullet (\maC^\infty _c (\maG))$ induced by an essential equivalence
$\varphi : \maK \to \maG$.
Hence in Section 2.3, we recall the results of Crainic \cite{Crainic;EtaleGpoid},
which computes the periodic homology of the convolution algebra in terms of the groupoid homology of the inertia groupoid $\maZ$.
This in turn enables us to apply Crainic and Moerdijk \cite{Moerdijk;EtaleHomology},
which describes the isomorphism of groupoid homology
$H _\bullet (\maK , \CC) \cong H _\bullet (\maG , \CC)$ in terms of a double complex constructed from an essential equivalence,
that we shall recall in section 3.3.
The isomorphism of periodic cyclic homology thus follows by combining \cite{Moerdijk;EtaleHomology},
with the well known fact that an essential equivalence $\varphi : \maK \to \maG$ induces an essential equivalence $\bar \varphi$,
of their respective inertia groupoids, which we give a proof in Section 3.2.
With these preliminaries,
in Section 4.1 we verify that $\varepsilon $ and $(f_y)_!$ (which is the composition of Thom isomorphism and exponential map)
commutes with the appropriate version of isomorphism induced by Morita equivalence,
hence proving Theorem \ref{CD6}.
Then we define our pairing in Section 4.2 and its invariance under Morita equivalence follows easily from Theorem \ref{CD6}.

Lastly, in Section 5, we return to the group action case.
We compute the index pairing between a function that is 1 on a conjugacy class and 0 outside,
with the Chern character of a Dirac type operator.
When all the non-trivial fixed points are isolated, we relate our index pairing
to the character values of the finite-group equivariant index through the Atiyah--Bott fixed-point formula.


\subsection*{Acknowledgements}
We would like to thank Raphael Ponge, Xiang Tang, and Bai-Ling Wang for inspirations and discussions.
HW is supported by the grants 23JC1401900 and NSFC 12271165,
and in part by the Science and Technology Commission of Shanghai Municipality (No.~22DZ2229014).

\subsection*{AI statement}
The ideas in this paper predate the major advances of the uses of AI in mathematics in the summer of 2026.
All proofs in this paper are solely developed by the authors.
The authors also confirm that all arguments, calculations,
and conclusions presented here were written and checked by the authors themselves.

\section{Cyclic homology of the convolution algebra of an (\'etale) Lie groupoid}

\subsection{Groupoids and their convolution algebras}
We review the definition of (\'etale and Lie) groupoids \cite{Ruan;Book, LN01, MackenzieBook05, MoerdijkFolBook, RenBook80, Sims20}.
A groupoid $\maG = \maG ^{(1)} \rightrightarrows \maG ^{(0)}$
consists of a set of objects $\maG ^{(0)}$ and a set of morphisms $\maG ^{(1)}$, endowed with structure maps:
\begin{enumerate}
\item
unit map $u : \maG^{(0)} \to \maG^{(1)}$;
\item
source and target maps $\bs , \bt : \maG ^{(1)} \to \maG ^{(0)}$;
\item
multiplication (or composition) $\mu : \maG ^{(2)}
:= \{ (\gamma _1 , \gamma _2 ) \in \maG ^{(1)} \times \maG ^{(1)} : \bs (\gamma _1 ) = \bt (\gamma _2) \} \to \maG ^{(1)}$
\item
inverse $\iota : \maG ^{(1)} \to \maG ^{(1)}$.
\end{enumerate}


These structural maps $\bs, \bt, \mu, u, \iota$ are required to satisfy the following properties
so that $\maG$ is a category:
\begin{enumerate}
\item
$\bs(\gamma _1 \gamma _2) = \bs(\gamma _2)$ and $\bt(\gamma _1 \gamma _2) = \bt(\gamma _1)$;
\item $(\gamma _1 \gamma _2) \gamma _3 = \gamma_1 (\gamma _2 \gamma _3)$, whenever they can be composed;
\item $\bs (u(x)) = x = \bt(u(x))$, $u(\bt(\gamma)) \gamma = \gamma = \gamma u(\bs(\gamma))$,
$\gamma \gamma ^{-1} = u (\bt(\gamma))$,  and $\gamma^{-1} \gamma = u(\bs(\gamma))$,
for any $\gamma \in \maG$ and $x \in \maG^{(0)}$.
\end{enumerate}
For $x \in \maG ^{(0)}$, we denote
$$ \maG _x^x := \{ \gamma \in \maG ^{(1)} : \bs (\gamma) = x = \bt (\gamma) \}, $$
it is a group, called the {\em isotropy group} at $x$.
However, if $U \subset \maG^{(0)}$, we simply write
$$\maG _U := \bs^{-1}(U) \cap \bt^{-1}(U),$$
which is a sub-groupoid, called the {\em reduction} of $\maG$ to $U$.

A {\em locally compact} groupoid is a groupoid $\maG$ equipped with a locally compact topology
such that the unit space $\maG^{(0)}$ is identified homeomorphically with a closed,
Hausdorff subset of $\maG$ via the unit map $u:\maG^{(0)}\rightarrow \maG$,
all structural maps $\bs, \bt, \mu, u, \iota$ are continuous,
and the source map $\bs$ is open.
In the present paper, all groupoids are Hausdorff.

A locally compact groupoid $\maG$ is called {\em \'etale} if $\bs$ and $\bt$ are both local homeomorphisms.
A {\em Lie} groupoid is a locally compact groupoid $\maG$ such that
$\maG$ and $\maG^{(0)}$ are Hausdorff manifolds,
the structural maps $\bs, \bt, \mu, u, \iota$ are smooth (in particular, $\iota$ is a diffeomorphism),
and $\bs, \bt: \maG \rightarrow \maG^{(0)}$ are both submersions.
A Lie or \'etale groupoid is called {\em proper} if the map $(\bs, \bt) : \maG \rightarrow \maG^{(0)} \times \maG^{(0)}$ is proper.
Most of the groupoids we consider in this paper are \'etale Lie groupoids,
which means that $\bs$ and $\bt$ are both local diffeomorphisms.

\begin{example}[Action groupoid]
Let $G$ be a discrete group acting (on the right) on a topological space $M$.
Set $\maG^{(1)}:= M \times G$ and $\maG^{(0)}:= M \cong M \times \{ e\}$, where
$e$ is the identity of the group $G$.
For all $x \in M$ and $g, h \in G$, define
$$\bs (x,g) = x g, \bt (x, g)= x, (x, g) (x g, h)=(x, g h), (x,g)^{-1}=(x g, g^{-1}).$$
Equipped with these structural maps,
$\maG^{(1)}\tto \maG^{(0)}$ becomes an \'etale groupoid.
We denote such an action groupoid by $M \rtimes G$.
\end{example}

To define a convolution algebra on a locally compact groupoid $\maG$,
we need a Haar system \cite{RenBook80}.
For an \'etale groupoid $\maG$ each $\bs ^{-1} (x)$ is discrete, and the counting measure is a Haar system.

\begin{dfn}
\label{Star}
Let $\maG$ be an \'etale Lie groupoid,
the convolution algebra $\maC ^\infty _c (\maG)$ is defined to be $C_c ^\infty (\maG ^{(1)})$ as a vector space,
but endowed the convolution product
$$ f_1 \star f_2 (\gamma) := \sum _{\bs (\gamma _1 ) = \bs (\gamma)} f_1 (\gamma \gamma _1 ^{-1}) f_2 (\gamma _1),$$
and involution
$$f ^* (\gamma) := \overline {f (\gamma ^{-1})} .$$
It is well-known that $\maC ^\infty _c (\maG)$ is an associative $*$-algebra with the above operations.
\end{dfn}

\begin{nota}
We shall denote the $*$-algebra of smooth,
compactly supported functions on a manifold $M$ with pointwise multiplication by $C ^\infty _c (M)$.
In particular when $\maG = \maG ^{(1)} \rightrightarrows \maG ^{(0)}$ is a Lie groupoid we write
$C ^\infty _c (\maG) := C^\infty _c (\maG ^{(1)})$ which differs from $\maC ^\infty _c (\maG)$ as an algebra.
\end{nota}

%

\begin{nota}
Let $\maG$ be a groupoid. We denote
$$\maG ^{(n)} := \{(\gamma _1 , \cdots, \gamma_n ) \in \maG^{(1)} \times \cdots \times \maG^{(1)} :
\bs(\gamma_i) = \bt(\gamma _{i+1}), i = 1, \cdots, n-1 \}.$$
\end{nota}

\subsection{Pre-cyclic objects and their cyclic homology}
We begin with recalling the basic definition of Hochschild and cyclic homology.
We mainly follow the standard reference \cite{Loday;Book}.

\begin{dfn}
\label{CyclicObj}
A pre-cyclic object is a simplicial object $\{ X _n \} _{n \geq 0} $, endowed with face maps
$$ d_i : X_n \to X _{n-1}, \quad i = 0, 1, \cdots, n, $$
and
$$t_{n+1} : X_n \to X_n,$$
satisfying the axioms:
\begin{enumerate}
\item
$d_i d_j = d _{j-1} d_i, \quad i < j $;
\item
$d _i t _{n+1} = \big\{
\begin{array}{ll}
t _n d _{i-1}, & 1 \leq i \leq n, \\
d _n, & i = 0 ;
\end{array}
$
\item
$t_{n+1}^{n+1} = 1$.
\end{enumerate}
Note that the last axiom $t_{n+1}^{n+1} = 1$ is equivalent to an action of $\ZZ _{n+1}$ on $X_n$,
defined on its canonical generator by $(-1)^n t_{n+1}$.
\end{dfn}

\begin{ex}
Let $A$ be an algebra over $\CC$ that is not necessarily unital.
On $A^\natural _n := \otimes ^{n+1} A$, define the maps
$ d_i : \otimes ^{n+1} A \to \otimes ^{n} A, i = 0, 1, \cdots, n $ and $t_{n+1} : \otimes ^{n+1} A \to \otimes ^{n+1} A$
by
\begin{align*}
d_i (a_0 \otimes \cdots \otimes a_n)
:=& a_0 \otimes \cdots \otimes a_{i-1} \otimes a_i a_{i+1} \otimes \cdots \otimes a_n, & 0 \leq i \leq n-1, \\
d_n (a_0 \otimes \cdots \otimes a_n)
:=& a_n a_0 \otimes \cdots \otimes a _{n-1}, & \\
t _{n+1} (a_0 \otimes \cdots \otimes a_n)
:=& a_n \otimes a_0 \otimes \cdots \otimes a_{n-1}.
\end{align*}
Then $(A^\natural _\bullet, d_i, t _{n+1}) $ is a pre-cyclic object.
If $A$ is unital one may define
$$s_i (a_0 \otimes \cdots \otimes a_n)
:= a_0 \otimes \cdots \otimes a_i \otimes 1 \otimes a_{i+1} \otimes \cdots \otimes a_n, \quad 0 \leq i \leq n.$$
Then $(A^\natural _n , d_i , s_i , t_{n+1})$ becomes a cyclic object.
\end{ex}

\begin{rem}
Our $t _{n+1}$ differs from Loday \cite[2.1.0]{Loday;Book} by a sign of $(-1)^n$.
\end{rem}

\begin{ex}
\label{CommutativeCyclic}
We shall consider topological algebras and take topological tensor products.
In particular, let $A = C _c ^\infty (M)$ for some manifold $M$.
We take completed tensor product (see\cite{HN77} and also \cite[Appendix B]{Tang;ProperLieCyclic}) so that
$$ C ^\infty _c (M) ^{\otimes (n+1)} := C _c ^\infty (M ^{n+1}) .$$
Then the pre-cyclic structure
\begin{align*}
d _i f (\gamma _0 , \cdots , \gamma _{n-1}) :=&
f (\gamma _0 , \cdots , \gamma_{i-1}, \gamma _i , \gamma _i , \gamma _{i+1}, \cdots , \gamma _{n-1}),
\quad 0 \leqslant i \leqslant n-1, \\
d _n f (\gamma _0 , \cdots , \gamma _{n-1}) :=&
f (\gamma _0 , \gamma_1 , \cdots , \gamma _{n-1}, \gamma_0 ), \\
t _{n+1} f (\gamma _0 , \cdots , \gamma _n) := & f (\gamma _1 , \cdots , \gamma _{n} , \gamma _0),
\end{align*}
extends the pre-cyclic structure on the algebraic tensor product.
\end{ex}

\begin{ex}
\label{CompletedTensor}
\cite[Equation (5)]{Nistor;EtaleCyclicHomology}
We extend the construction of Example \ref{CommutativeCyclic} to the smooth compactly supported convolution algebra,
$\maC^{\infty}_c (\maG)$,
of an \'etale groupoid $\maG$.
Again we consider the topological tensor product
$$ \maC^{\infty}_c (\maG) ^{\otimes (n+1)} := \maC^{\infty}_c (\maG ^{n+1}) .$$
Then $\maC^{\infty}_c (\maG) ^{\otimes (n+1)}$ is endowed with the pre-cyclic structure
\begin{align*}
d _i f (\gamma _0 , \cdots , \gamma _{n-1}) :=&
\sum _{\bt (\gamma) = \bt (\gamma _i)}
f (\gamma _0 , \cdots , \gamma_{i-1}, \gamma , \gamma ^{-1} \gamma _i , \gamma _{i+1}, \cdots , \gamma _{n-1}),
0 \leqslant i < n; \\
d _n f (\gamma _0 , \cdots , \gamma _{n-1}) :=&
\sum _{\bt (\gamma) = \bt (\gamma _0)}
f (\gamma ^{-1} \gamma _0 , \gamma _1 , \cdots , \gamma _{n-1}, \gamma ); \\
t _{n+1} f (\gamma _0 , \cdots , \gamma _n) := & f (\gamma _1 , \cdots , \gamma _{n} , \gamma _0).
\end{align*}
We denote such pre-cyclic object by $\maC^{\infty}_c (\maG) ^\natural$.
\end{ex}

Given a pre-cyclic object $(X_n , d_i , t_{n+1})$, define the maps
$$ b, b' : X_n \to X_{n-1} , \quad b' := \sum _{j=0}^{n-1} (-1)^j d _j, \quad b := b' + (-1)^n d _n, $$
$$ \epsilon := 1 - (-1)^n t _{n+1}, \quad N := \sum _{j=0}^n (-1)^{j n} t_{n+1} ^j .$$

\begin{dfn}
The Hochschild homology of a pre-cyclic object $(X_n , d_i , t_{n+1})$ is defined to be the total homology of the bi-complex
$$
\begin{CD}
@VVV @VVV \\
X _{n+1} @<{\epsilon}<< X _{n+1}  \\
@V{b}VV @V{-b'}VV  \\
X _{n} @<{\epsilon}<< X _{n} \\
\vdots @. \vdots  \\
@VVV @VVV \\
X _{0} @<{\epsilon}<< X _{0}
\end{CD}
$$
and we denote the Hochschild homology of $(X_n , d_i , t_{n+1})$ by $H H _k (X)$.
\end{dfn}

Recall that for $(A^\natural _\bullet, d_i, t _{n+1}) $,
if $A$ is H-unital, then the $b'$ column is acyclic and the Hochschild homology is just the homology of the complex
$ (A^\natural _\bullet , b)$.

\begin{dfn}
\label{PreCyclicHom}
The cyclic homology of a pre-cyclic object $(X_n , d_i , t_{n+1})$ is defined to be the total homology of the Connes-Tsygan bi-complex
\begin{equation}
\label{CyclicComplex}
\begin{CD}
@VVV @VVV @VVV @. \\
X _{n+1} @<{\epsilon}<< X _{n+1} @<{N}<< X _{n+1} @<{\epsilon}<< \\
@V{b}VV @V{-b'}VV @V{b}VV @. \\
X _{n} @<{\epsilon}<< X _{n} @<{N}<< X _{n} @<{\epsilon}<< \\
\vdots @. \vdots @. \vdots @. \\
@VVV @VVV @VVV @. \\
X _{0} @<{\epsilon}<< X _{0} @<{N}<< X _{0} @<{\epsilon}<<
\end{CD}
\end{equation}
and we denote the cyclic homology of $(X_n , d_i , t_{n+1})$ by $H C _k (X)$.
Explicitly the cyclic homology is the homology of the complex
$$ ( X_n \oplus \cdots \oplus X_0, (b + \epsilon) \oplus (-b'+N) \oplus \cdots ). $$
\end{dfn}

\begin{dfn}
The periodic cyclic homology $H P _k (X)$ is defined similarly.
Instead of \eqref{CyclicComplex} consider the bi-complex \cite[Section 5.1.1]{Loday;Book}
$$
\begin{CD}
@. @VVV @VVV @VVV @. \\
@<<< X _{n+1} @<{\epsilon}<< X _{n+1} @<{N}<< X _{n+1} @<<< \\
@. @V{b}VV @V{-b'}VV @V{b}VV @. \\
@<<< X _{n} @<{\epsilon}<< X _{n} @<{N}<< X _{n} @<<< \\
@. \vdots @. \vdots @. \vdots @. \\
@. @VVV @VVV @VVV @. \\
@<<< \underbrace{X _{0}}_{\deg 0} @<{\epsilon}<< \underbrace{X _{0}}_{\deg 1} @<{N}<< X _{0} @<<<
\end{CD}
$$
which induces the period-2 total complex
$$ \leftarrow \underbrace{\prod _{k \geq 0} X_k}_{\deg 0} \xleftarrow{((b + \epsilon) \oplus (-b'+N) \oplus \cdots)}
\underbrace{\prod _{k \geq 0} X_k}_{\deg 1} \xleftarrow{((-b'+N) \oplus (b + \epsilon) \oplus \cdots)}
\prod _{k \geq 0} X_k \leftarrow ,$$
where the infinite Cartesian product $\prod _{k \geq 0} X_k $ is the space of sequences $\{ x_k \}_{k \geq 0}$ with $x_k \in X_k$.
Its homology is called the periodic cyclic homology of $(X_n , d_i , t_{n+1})$, denoted by $H P _k (X)$ ($k = 0, 1$).
\end{dfn}

Dually, the cyclic cohomology $H C ^k (X)$ and periodic cyclic cohomology $H P ^k (X)$ are defined by applying
$\operatorname{Hom}(-,\CC)$ to the corresponding cyclic and periodic cyclic complexes; see \cite{Loday;Book}.

\begin{rem}
For a cyclic object $(X _n , d_i , s_i , t_{n+1})$ the cyclic homology is equivalently defined by the $b$-$B$ bi-complex,
where $b$ is the same as Definition \ref{PreCyclicHom} and
$$B := \varepsilon s _0 N .$$
For the general case one can also adjoint a unit and consider the $b$-$B$ complex.
\end{rem}

\subsection{Groupoid homology and cyclic homology of convolution algebra}
\label{Crainic}
Let $\maG$ be a Hausdorff \'etale groupoid.
In this section, we recall how one computes the (periodic) cyclic homology $H C _\bullet (\maC^\infty_c (\maG))$
in terms of groupoid homology,
following Crainic \cite{Crainic;EtaleGpoid}.
More details can be found in Appendix A.

\begin{dfn}
Let $\maG$ be an \'etale groupoid, $\maA$ be a c-soft $\maG$-sheaf of complex algebras.
The algebraic cross product is $ \maA \rtimes \maG := \Gamma _c (\maG, \bs^* \maA) $,
endowed with convolution product (determined by the germs)
$$ (u \star v) _\gamma := \sum _{\bs (\gamma') = \bs (\gamma)} (u _{\gamma (\gamma')^{-1}} \cdot \gamma') v _{\gamma'} .$$
\end{dfn}

\begin{ex}
\label{CrossProd}
When $\maA = C ^\infty_{\maG ^{(0)}}$, the sheaf of germs of smooth functions on $\maG ^{(0)}$ (see Example \ref{SmoothGerms}),
$\maA \rtimes \maG $ is just $\maC^\infty_c (\maG)$ as an algebra.

In this case we let $\maA ^{\otimes n}$ ($= (C ^\infty _{\maG ^{(0)}}) ^{\otimes n}$) be the sheaf over $(\maG ^{(0)}) ^n$
($n$-fold Cartesian product), as in Remark \ref{CommutativeCyclic},
where $\maA ^{\otimes n} (U)$ is the algebra of smooth functions on $U ^n$.
\end{ex}


\begin{dfn}
Given an \'etale groupoid $\maG$, Burghelea's space is defined to be
$$ B ^{(n)} := \{ (\gamma _0, \cdots , \gamma _n) \in \maG ^{(n+1)} : \bt (\gamma _0) = \bs (\gamma _n) \}, \quad n = 0, 1, \cdots .$$
\end{dfn}

%
%
%
%

\begin{dfn}\label{inertialgroupoid}
The groupoid $\maG$ acts on $B ^{(0)}$ from the right by conjugation
$$ \mu \cdot \gamma := \gamma ^{-1} \mu \gamma $$
for $\mu \in B ^{(0)}, \gamma \in \maG ^{(1)}$ satisfying $\bs (\mu) = \bt (\gamma )$.

The inertia groupoid is defined to be the action groupoid $\maZ := B ^{(0)} \rtimes \maG$.
\end{dfn}

Observe that because $\maG ^{(0)}$ is open and closed in $B ^{(0)}$,
we may write $B ^{(0)}$ into a disjoint union of two open and closed invariant subsets, namely
\begin{align*}
B ^{(0)} _\ZZ := & \{ \gamma : \gamma ^ r \in \maG ^{(0)} \text { for some $r \in \ZZ \setminus \{ 0 \}$} \} \\
B ^{(0)} _\infty := & \{ \gamma : \gamma ^ r \not \in \maG ^{(0)} \text { for all $r \in \ZZ \setminus \{ 0 \}$} \},
\end{align*}
elements of $B ^{(0)} _\ZZ $ are called elliptic, while elements of $B ^{(0)}_\infty $ are called hyperbolic.
Correspondingly, we denote
$$ \maZ _\ZZ := B ^{(0)} _\ZZ \rtimes \maG , \quad \maZ _\infty := B ^{(0)} _\infty \rtimes \maG .$$

There is a map $\theta : B ^{(0)} \to \maZ ^{(1)} $
$$ \theta (z ) := (z , z ) \in B ^{(0)} \times _\bt \maG ^{(1)} .$$
This defines a cyclic groupoid structure on $\maZ$.
Let $\maN$ be the quotient of $\maZ $ by the relation
$$ (z , z ^n \gamma ) \sim (z , \gamma ), \quad z \in B ^{(0)}, \gamma \in \maG ^{(1)}, n \in \ZZ , $$
and similarly $\maN _\ZZ $ and $\maN _\infty$ for $\maZ _\ZZ$ and $\maZ _\infty$ respectively,
$\maN$ is still a groupoid.

\begin{dfn}
\label{TwistedLoop}
Let $\maA = C ^\infty _{\maG ^{(0)}}$,
we pull back to $B ^{(0)}$ by the restriction of the source map $\bs _\maG$ on $\maG$ to $\maB ^{(0)}$,
i.e., let
$$\maA _{B ^{(0)}} := \bs _{\maG} \vert _{\maB ^{(0)}} ^* \maA .$$
(see \cite[3.24]{Crainic;EtaleGpoid}.)
Define $\maC (\maZ , \theta , \maA)$ to be $\Gamma _c \big(\maZ , \bs ^* \maA _{\maB ^{(0)}} ^{\otimes (\bullet+1)} \big)$,
where $\bs$ here is the source map of $\maZ$.
Endow $\maC (\maZ , \theta , \maA)$ with the pre-cyclic structure
\begin{align*}
\tilde d_0 f _{(\gamma _1, \cdots , \gamma _n)}
:=& d_0 f _{(\gamma _2 , \cdots , \gamma _n)} \cdot \gamma _1,\\
\tilde d_i f _{(\gamma _1, \cdots , \gamma _n)}
:=& d_i f _{(\gamma _1 , \cdots , \gamma _i \gamma _{i+1} , \cdots , \gamma _n)} , \quad 1 \leq i < n, \\
\tilde d_n f _{(\gamma _1, \cdots , \gamma _n)}
:=& d_n f _{(\gamma _1 , \cdots , \gamma _{n-1})}, \\
\tilde t f _{(\gamma _1 , \cdots , \gamma _n)}
:=& (t f _{((\gamma _1 \cdots \gamma _n)^{-1} \theta , \gamma _1 , \cdots , \gamma _{n-1})})
\cdot \theta ^{-1} \gamma _1 \cdots \gamma _n,
\end{align*}
where $d_i$ and $t$ are defined as in Remark \ref{CommutativeCyclic} (on each stalk).
\end{dfn}

Since $\maC (\maZ , \theta , \maA)$ is a pre-cyclic object we shall consider its Hochschild, cyclic and periodic cyclic homology.

\begin{lem}
\label{Reduction}
\cite[Proposition 3.36, 4.1]{Crainic;EtaleGpoid}
We have canonical isomorphisms
\begin{align*}
H H _\bullet (\maC (\maZ , \theta , \maA)) \cong & H H _\bullet (\maC ^\infty _c (\maG)), \\
H C _\bullet (\maC (\maZ , \theta , \maA)) \cong & H C _\bullet (\maC ^\infty _c (\maG)), \\
H P _\bullet (\maC (\maZ , \theta , \maA)) \cong & H P _\bullet (\maC ^\infty _c (\maG)) .
\end{align*}
\end{lem}

\begin{proof}[Sketch of proof of Lemma \ref{Reduction}]
Here we briefly recall the construction of the isomorphism above.
In the following we shall write $\maA := C^\infty _{\maG ^{(0)}}$ as in Example \ref{CrossProd}.

To begin with, on $\maG$ define $ \sigma _n : B ^{(n)} \to (\maG ^{(0)})^{n+1}$ to be
$$ \sigma _n (\gamma _0, \cdots , \gamma _n) := (\bs _\maG (\gamma _0), \cdots, \bs _\maG (\gamma _n)),$$
and let
\begin{equation}
\label{AlgCross}
(\maA \rtimes \maG) ^\natural (n) := \Gamma _c (B ^{(n)} , \sigma _n ^* (\maA ^{\otimes n+1})) .
\end{equation}
By \cite[Lemma 3.30]{Crainic;EtaleGpoid}, $(\maA \rtimes \maG) ^\natural$ is endowed with a (pre-)cyclic structure.

By \cite[Lemma 3.34]{Crainic;EtaleGpoid} and \cite[Section 3]{Nistor;EtaleCyclicHomology},
given $[c] \in H P _\bullet (\maC _c ^\infty (\maG))$,
restricting $c$ to its sheaf over $B ^{(n)}$ induces a quasi-isomorphism from
$\maC^{\infty}_c (\maG ) ^\natural$ to $(\maA \rtimes \maG) ^\natural$.

Next, one identifies $(\maA \rtimes \maG) ^\natural$ with $\maC (\maZ , \theta , \maA)$
by the bijection \cite[Proposition 3.36]{Crainic;EtaleGpoid}
\begin{align}
\label{Twist}
(z_0, \cdots, z_n &; \gamma _0, \cdots, \gamma _n) \\ \nonumber
\mapsto & (z_0 \gamma_1 \cdots \gamma_n \gamma_0, z_0 \gamma_2 \cdots \gamma_n \gamma_0, \cdots, z_n \gamma _0;
\gamma_1 \cdots \gamma_n \gamma_0, \gamma_1, \cdots, \gamma_n). \qedhere
\end{align}
\end{proof}

Then one proceeds by splitting $\maZ$ into its elliptic and hyperbolic parts $\maZ _\ZZ$ and $\maZ _\infty$.
For the elliptic part, one has
$$ H H _\bullet (\maC (\maZ _\ZZ , \theta , \maA)) \cong H H _\bullet (\maC (\maN _\ZZ , 1 , \maA)). $$
Then by the Eilenberg-Zilber theorem cited in \cite{Crainic;EtaleGpoid},
$H H _\bullet (\maC (\maN _\ZZ , 1 , \maA)) $ is computed by the groupoid (hyper-)homology
\cite[Proposition 3.26]{Crainic;EtaleGpoid}:
$$ H H _\bullet (\maC (\maN _\ZZ , 1 , \maA))
\cong H _\bullet \big(\maN _\ZZ , ((\maA _{B ^{(0)}_\ZZ} )^\natural _\bullet , b) \big).$$
Then the Hochschild-Kostant-Rosenberg theorem can be applied to the pre-cyclic object coming from the commutative algebra,
$(\maA _{B ^{(0)}_\ZZ} )^\natural$, one has quasi-isomorphism with the de Rham complex \cite[Lemma 3.5]{Crainic;EtaleGpoid}
$$ \big((\maA _{B ^{(0)}_\ZZ} )^\natural _\bullet , b \big) \to \big( \Omega _{B ^{(0)}_\ZZ} ^\bullet , 0 \big),$$
that furthermore makes the $B$-complex (which is well defined at each stalk because each $\maA _z$ is unital)
compatible with the de Rham differential.
In particular we have quasi-isomorphism between the periodic bi-complexes of sheaves:
$$
\begin{CD}
@. @VbVV @VbVV @VbVV \\
\cdots @<B<< (\maA _{B ^{(0)}} )^\natural_2 @<B<< (\maA _{B ^{(0)}} )^\natural_1 @<B<< (\maA _{B ^{(0)}} )^\natural_0 \\
@. @VbVV @VbVV \\
\cdots @<B<< (\maA _{B ^{(0)}} )^\natural_1 @<B<< (\maA _{B ^{(0)}} )^\natural_0 \\
@. @VbVV \\
\cdots @<B<< (\maA _{B ^{(0)}} )^\natural_0
\end{CD}
$$
and
$$
\begin{CD}
@. @V{0}VV @V{0}VV @V{0}VV \\
\cdots @<{d _{\DR}}<< \Omega ^2 _{B^{(0)}} @<{d _{\DR}}<< \Omega ^1 _{B^{(0)}} @<{d _{\DR}}<< \Omega ^0 _{B^{(0)}} \\
@. @V{0}VV @V{0}VV \\
\cdots @<{d _{\DR}}<< \Omega ^1 _{B^{(0)}} @<{d _{\DR}}<< \Omega ^0 _{B^{(0)}} \\
@. @V{0}VV \\
\cdots @<{d _{\DR}}<< \Omega ^0 _{B^{(0)}}
\end{CD}
.
$$
Hence the periodic cyclic homology $ H P _\bullet (\maC (\maZ _\ZZ , \theta , \maA)) $ is isomorphic to the homology of
$$ \Big( \prod _{p} B _p \big(\maG , \bigoplus _k \Omega ^{\bullet + p + 2k} _{B ^{(0)}} \big), \partial + d _{\DR} \Big). $$
Since
$$ 0 \to \CC \to \Omega^0 _{B ^{(0)}} \xrightarrow{d_{\DR}} \Omega ^1 _{B^{(0)}} \to \cdots $$
is an exact sequence of c-soft sheaves,  we conclude that:
\begin{lem}
\label{Elliptic}
\cite[Theorem 4.4]{Crainic;EtaleGpoid}
The periodic cyclic homology
$ H P _\bullet (\maC (\maZ _\ZZ , \theta , \maA)) $ is isomorphic to the groupoid homology
$$ \prod _k H _{\bullet + 2 k} (\maN _\ZZ , \CC) \cong \prod _{k} H _{\bullet + 2 k} (\maZ _\ZZ , \CC).$$
\end{lem}

As for the hyperbolic case, we have
\begin{lem}
\cite[Corollary 4.7]{Crainic;EtaleGpoid}
\label{FiniteDim}
If $\maN _\infty$ has finite cohomological dimension, then
$$ H P _\bullet (\maC (\maZ _\infty , \theta , \maA)) \cong 0 .$$
\end{lem}

\begin{lem}
\cite[Proposition 4.6]{Crainic;EtaleGpoid}
If every $z \in B ^{(0)} _\infty $ is stable, then
$$ H C _\bullet (\maC (\maZ _\infty, \theta , \maA))
\cong H _{p-q=\bullet} \Big(\maN _\infty , (\Omega _{B ^{(0)}_\infty} ^\bullet , 0, d _{\DR} + \cap _e) \Big),$$
where $\cap _e$ denotes taking cap product with the Euler class.
\end{lem}

\subsection{Decomposition of cyclic homology of action groupoids}
We specialize the construction in the previous section to action groupoids.
Let $G$ be a discrete group acting on a manifold $M$ from the right.
Let $\maG := M \rtimes G \rightrightarrows M$ be the action groupoid.

Note that there is an equivalent way to study the convolution algebra, namely one has
$$ \maC^\infty_c (\maG) \cong (C^\infty _c (M)) \rtimes G $$
as algebras. Also one has isomorphism of the bar complexes (Definition \ref{BarComplex})
\begin{equation}
\label{ActionBar}
B _n (\maG , C^\infty _M) \cong B _n (G, \varrho _! C^\infty _M) ,
\end{equation}
where $C^\infty_M$ is the sheaf of germs of smooth functions on $M$,
and $\varrho$ is the (constant) map from $M$ to the single point that is the unit of $G$.

The loop space can be explicitly written down
\begin{equation}
\label{Loop}
B ^{(0)} = \{ (x, l) \in M \times G : x \cdot l = x \} ,
\end{equation}
and the inertia groupoid is just
$$ \maZ = B ^{(0)} \rtimes \maG = B ^{(0)} \rtimes G.$$
Other than just the elliptic and hyperbolic parts,
$B ^{(0)}$ can be further decomposed in this case, according to the conjugacy class of $G$.
For each $l \in G$, let $M^l$ be the sub-manifold fixed by $l$.
Then for each conjugacy class $\langle g \rangle $
\begin{equation}
B ^{(0)} _{\langle g \rangle} := \{ (x, l) \in B ^{(0)} : l \in \langle g \rangle \}
\end{equation}
is a $G$-invariant, open and closed subset.
Moreover $B^{(0)} _{\langle g \rangle} \rtimes G$ is Morita equivalent to $M^g \rtimes Z_g $
(here $Z _g \subset G$ denotes the centralizer subgroup of $g$).

Observe that $B^{(0)} _{\langle g \rangle} \subset B^{(0)} _\ZZ$ if $g$ has finite order,
and $B^{(0)} _{\langle g \rangle} \subset B^{(0)} _\infty$ if $g$ has infinite order.

For the elliptic case, we have
$$ H P _\bullet (\maC (\maZ _\ZZ , \theta , \maA))
= \bigoplus _{\langle l \rangle} H P _\bullet (\maC (\maZ _{\langle l \rangle} , \theta , \maA)) .$$
Then Lemma \ref{Elliptic} can be applied to each component.
\begin{cor}
\label{HPSplitting}
We have
\begin{align*}
H P _\bullet (\maC (\maZ _{\langle l \rangle} , \theta , \maA))
\cong & \prod _k H _{\bullet + 2 k} (B^{(0)}_{\langle l \rangle} \rtimes G , \CC) \\
\cong & \prod _k H _{\bullet + 2 k}
(B^{(0)}_{\langle l \rangle} \rtimes G , (\Omega _{B^{(0)}_{\langle l \rangle}}^\bullet, d _{\DR}))\\
\cong & \prod _k H _{\bullet + 2 k} (G , (\varrho _! \Omega _{B^{(0)} _{\langle l \rangle}} ^\bullet , d _{\DR}))
\end{align*}
(see Definition \ref{ChainCyclic} for the definition of
$H _{\bullet + 2 k} (G , (\varrho _! \Omega _{B^{(0)} _{\langle l \rangle}} ^\bullet , d _{\DR}))$).
\end{cor}
\begin{proof}
The first isomorphism is Lemma \ref{Elliptic}.
The second isomorphism follows from the fact that
$\Omega _{B^{(0)} _{\langle l \rangle}} ^\bullet , d _{\DR})$ is a resolution of
$\CC$ as a $B^{(0)} _{\langle l \rangle} \rtimes G$-sheaf.
The last isomorphism follows from \eqref{ActionBar}.
\end{proof}

\subsection{Thom isomorphism for action groupoids}
Let $M$ be an oriented manifold, $E \xrightarrow{\wp} M$ be an oriented real vector bundle of rank $r$.
Then the Thom isomorphism
$$ \omega \mapsto \tau \wedge \wp ^* \omega , $$
where $\tau \in \Omega ^r (E)$ is any representative of the Thom class,
induces a quasi-isomorphism from the compactly supported de Rham cohomology $H ^\bullet _c ( M ) $ to $H ^{\bullet + r} _c ( E )$.
The inverse isomorphism is given by fiber-wise integration.

Let $G$ be a discrete group acting on $M$ and $E$ by orientation-preserving diffeomorphisms,
so that the action commutes with $\wp$.
We set
$$\maG := M \rtimes G, \quad \tilde{\maG} := E \rtimes G = E \rtimes \maG.$$
We shall compare the periodic cyclic homology of $\maC^\infty _c (\maG)$ and $\maC^\infty _c (\tilde \maG)$.
The loop space of $\maG$  is given by \eqref{Loop},
and for $\tilde \maG$ the loop space is similarly given by
$$ \tilde B ^{(0)} = \{ (v, l) \in E \times G : v \cdot l = v \} ,$$
and splits into a disjoint union $\tilde B ^{(0)} _{\langle g \rangle} $ by conjugacy classes.
Moreover, observe that the vector bundle projection $\wp$,
acting on the $E$ component in $\tilde B$, defines a vector bundle structure
$\tilde B ^{(0)} _{\langle g \rangle} \to B^{(0)} _{\langle g \rangle}$
(and also $E ^g \to M^g$) with possibly different rank on each connected component.
For the Thom isomorphisms below, assume that each sector bundle
$\tilde B ^{(0)} _{\langle l \rangle} \to B^{(0)} _{\langle l \rangle}$
for finite-order $l \in G$ is $G$-equivariantly oriented.
Without this assumption, the Thom isomorphism involves the orientation local system.

Suppose that $l \in G$ has finite order.
Then for each $x \in M^l$, $l$ maps $E _x $ to itself,
and is given by a real matrix of positive determinant (because the action is orientation preserving).
Since $l$ has finite order, its nonreal eigenvalues occur in conjugate pairs and its real eigenvalues are $1$ or $-1$.
The positive determinant implies that the $-1$-eigenspace has even dimension.
Thus the sum of the dimensions of the eigenspaces with eigenvalue different from $1$ is even, and
$$\dim (E ^l)_x \equiv \rk (E) \pmod 2.$$

Because the compactly supported de Rham complex computes the groupoid homology,
the Thom isomorphism for the sector bundle induces a quasi-isomorphism, so that we have
$$ H _{\bullet + 2 k} (B^{(0)}_{\langle l \rangle} \rtimes G , (\Omega _{B^{(0)}_{\langle l \rangle}}^\bullet, d _{\DR}))
\cong H _{\bullet - r_l + 2 k}
(\tilde B^{(0)}_{\langle l \rangle} \rtimes G , (\Omega _{\tilde B^{(0)}_{\langle l \rangle}}^\bullet, d _{\DR})),$$
where $r_l$ is the rank of $E ^l$ over $M^l$, understood componentwise if this rank varies
(see Example \ref{ThomEx} for more about the contractible case).

Assuming that $H P _\bullet (\maC (\maZ _\infty , \theta , \maA)) \cong 0$ (hence leaving only the elliptic part),
it follows from Corollary \ref{HPSplitting} that
\begin{align}
\label{Thom}
H P _{\bullet} (C ^\infty _c (M) \rtimes G) \cong & H P _\bullet (\maC^\infty _c (\maG)) \\ \nonumber
\cong & H P _{\bullet - \rk (E)} (\maC^\infty _c (\tilde \maG)) \cong H P _{\bullet - \rk (E)} (C^\infty _c (E) \rtimes G).
\end{align}

\begin{dfn}
We call \eqref{Thom} the Thom isomorphism and denote it by
$$ \maT _E : H P _{\bullet} (C^\infty _c (M) \rtimes G) \to  H P _{\bullet - \rk (E)} (C^\infty _c (E) \rtimes G). $$
\end{dfn}

\subsection{The pushforward map}
Let $G$ be a discrete group acting on $M, N$ from the right.
We assume that the loop spaces
$$\{ (x, l) \in M \times G : x \cdot l = x \} \text{ and } \{ (y, l) \in N \times G : y \cdot l = y \} $$
are both oriented, and that $G$ acts on these loop spaces (Definition \ref{inertialgroupoid}) preserving orientation.
Let $f : M \to N$ be a $G$-equivariant map.
Observe that $f$ maps $M^l$ to $N^l$, and we assume that the induced map on the loop spaces preserves orientation.
Recall that P. Carrillo Rouse, B. L. Wang, and H. Wang constructed a pushforward map in \cite{CRWW23}
$$f_! : \bigoplus H P _\bullet (\maC^\infty _c (M \rtimes G))
\to \bigoplus H P _{\bullet + \dim M + \dim N} (\maC^\infty _c (N \rtimes G)), $$
whose construction we recall here (note that we consider right action instead).

\begin{dfn}
Let $\maH = \maH ^{(1)} \rightrightarrows \maH ^{(0)}$ be a Lie groupoid,
$G$ be a discrete group acting on $\maH$ from the right by groupoid automorphisms.
Then we denote by $\maH \rtimes G$ the action groupoid.
Recall that as manifolds,
$$(\maH \rtimes G)^{(1)} := \maH ^{(1)} \times G, \quad (\maH \rtimes G)^{(0)} := \maH ^{(0)} $$
and the source, target, and multiplication maps are given by
$$ \bs _{\maH \rtimes G} (\gamma, g) := \bs _\maH (\gamma \cdot g), \quad
\bt _{\maH \rtimes G} (\gamma, g) := \bt _\maH (\gamma), \quad
(\gamma , g) (\gamma' , g') := (\gamma ( \gamma' \cdot g^{-1}), g g').$$
\end{dfn}

Define the vector bundle
$$E := \{ (v, w, x) : (v, x) \in T M, w \in T _{f (x)} N \} \xrightarrow{\wp} f ^* T N,$$
and $E'$ be its dual bundle (over $f^* T N$).
Here $\wp(v,w,x)=(w,x)$, so the fiber of $E$ over $(w,x)$ is $T_xM$ and the fiber of $E'$ is $T_x^*M$.
For the Thom map in item (1), regard the total space of $E'$ as a vector bundle over $M$ via the projection $f^*TN \to M$;
its fiber over $x$ is $T_{f(x)}N \oplus T_x^*M$.
Fiber-wise vector addition endows $E$ with a Lie groupoid structure with unit space $f ^* T N$.
Define $C^\infty _s (E)$ to be the space of smooth functions $\phi \in C^\infty (E)$
such that $\wp (\supp (\phi)) \subset f ^* T N$ is compact and $\phi$ is rapidly decaying (Schwartz) on each fiber.

The convolution algebra $\maC^\infty _s (E \rightrightarrows f^* T N)$ is defined to be $C^\infty _s (E)$ as a vector space,
but endowed with fiber-wise convolution product on $E \rightrightarrows f^* T N$.
Fiber-wise Fourier transform is an algebra isomorphism from $\maC^\infty _s (E \rightrightarrows f^* T N)$ to $C^\infty _s (E')$,
the latter endowed with point-wise multiplication.
It follows that we have
$$ H P _\bullet (C^\infty _s (E') \rtimes G)  \cong H P _\bullet ( \maC^\infty_s (E \rightrightarrows f^* T N) \rtimes G) .$$
Note that $\maC^\infty _s (E \rightrightarrows f^* T N) \rtimes G$ is the (Schwartz) convolution algebra of
$E \rtimes G \rightrightarrows f^* T N$.

On the other hand, one can define a groupoid
$$ \maD := E \sqcup ((M \times M \times N) \times (0, 1]) \rightrightarrows f^* T N \sqcup (M \times N \times (0, 1]),$$
that can be furthermore endowed with a Lie structure,
and we shall study $\maD \rtimes G$.
In particular we consider two invariant sub-groupoids of $\maD \rtimes G$, namely
$$ E \rtimes G \rightrightarrows f ^* T N
\text{ and } ((M \times M \times N) \times \{1\}) \rtimes G \rightrightarrows M \times N \times \{ 1 \} .$$
Restriction to these invariant sub-groupoids gives algebra homomorphisms $e _0 , e _1$ from
$\maC^\infty _s (\maD \rtimes G)$ (the convolution algebra that is Schwartz on $E \rtimes G$ and compactly supported on each
$M \times M \times \{ y \}, y \in N$) to
$$ \maC^\infty _s (E \rtimes G \rightrightarrows f ^* T N) \text{ and respectively }
 \maC^\infty _c ((M \times M \times N) \rtimes G) .$$
Hence one has induced maps
\begin{align*}
e_0 : H P _\bullet (\maC^\infty _s (\maD \rtimes G)) \to & H P _\bullet (\maC^\infty _s (E \rtimes G \rightrightarrows f ^* T N)) \\
e_1 : H P _\bullet (\maC^\infty _s (\maD \rtimes G)) \to & H P _\bullet (\maC^\infty _c ((M \times M \times N) \rtimes G))
\cong H P _\bullet (\maC^\infty _c (N \rtimes G)),
\end{align*}
where the last isomorphism is defined by Morita equivalence (see Section 3 below, or \cite{CRWW23} for more details).
Using homotopy invariance and six terms exact sequence for periodic homology,
it can be shown that $e_0$ is an isomorphism.
Indeed $(e_0)^{-1}$ can be explicitly written in terms of the exponential map (see below).

\begin{dfn}
\label{Pushforward}
The pushforward map $f_! : H P _\bullet (\maC^\infty _c (M \rtimes G)) \to H P _\bullet (\maC^\infty _c (N \rtimes G))$
is the composition of:
\begin{enumerate}
\item
The Thom isomorphism
$$\maT _{E'} : H P _\bullet (\maC^\infty _c (M) \rtimes G)
\to H P _{\bullet + \dim M + \dim N} (C^\infty _s (E') \rtimes G)),$$
where the total space of $E'$ is regarded as a vector bundle over $M$;
\item
Fourier transform $H P _\bullet (C^\infty _s (E') \rtimes G)
\cong H P _\bullet (\maC^\infty _s (E \rtimes G \rightrightarrows f ^* T N) \rtimes G)$;
\item
Morphism induced by deformation groupoid
$e_1 \circ (e_0)^{-1} : H P _\bullet (\maC^\infty _s (E) \rtimes G) \to H P _\bullet (\maC^\infty _c ((M \times M \times N) \rtimes G))$;
\item
Morita equivalence
$H P _\bullet (\maC^\infty _c ((M \times M \times N) \rtimes G)) \to H P _\bullet (\maC^\infty _c (N \rtimes G))$.
\end{enumerate}
\end{dfn}

\begin{rem}
\label{NonProper}
Although \cite{CRWW23} generally assumes proper action (which is natural in the context of Baum-Connes conjecture),
the construction of the pushforward map, as described above, does not require the $G$-action to be proper.
Also, the orientation assumption can be weakened,
for example when $f: M \to N$ is an open embedding, $E = T M \oplus T M$ which is always oriented,
and hence the Thom isomorphism is well defined.
\end{rem}

We recall some basic properties of the pushforward map in \cite{CRWW23}.

\begin{lem}
\label{PushLem1}
Let $M$, $N$, and $L$ be $G$-proper cocompact manifolds, and let
$f_1:M\to N$ and $f_2:N\to L$ be $G$-equivariant, $K$-oriented smooth maps.
Then the pushforward is functorial:
$$ (f_2)_! \circ (f_1)_! = (f_2 \circ f_1)_! .$$
\end{lem}
This is the periodic cyclic homology identity proved in
\cite[Theorem 4.1, Equation (4.7)]{CRWW23}.

\begin{lem}
Let $F\to M$ be a $G$-proper, $G$-equivariant $\operatorname{spin}^c$ vector bundle
over a $G$-proper cocompact manifold, and let $s_0:M\to F$ be its zero section.
Then the pushforward along $s_0$ is the Thom isomorphism:
$$ (s_0)_! = \maT_F .$$
Consequently, if $h:F\to N$ is a $G$-equivariant, $K$-oriented smooth map
between $G$-proper cocompact manifolds, then
$$ (h\circ s_0)_! = h_!\circ \maT_F .$$
\end{lem}
The first identity is the zero-section case in \cite[Section 3.1]{CRWW23},
and the second follows from \cite[Theorem 4.1]{CRWW23}.


We explicitly describe the case when $M = \{ x _0 \}$ is the single point space.
Let $\maG = \maG ^{(1)} \rightrightarrows \maG ^{(0)}$ be a Lie groupoid and
$U $ be any open subset of $\maG ^{(0)}$.
Recall that $\maG _U := \bs ^{-1} _\maG (U) \cap \bt ^{-1} _\maG (U)$,
which is an open subgroupoid of $\maG$.
We shall further assume that there is a fixed, $G$-invariant affine connection $\nabla$ on $T N$
(when the $G$-action is proper one may let $\nabla$ be the Riemannian connection of an invariant metric).
Let $y _0 : f (x _0) \in N$.
Note that $G$ fixes $y$.

Since $M=\{x_0\}$, the bundles $E\to T_{y_0}N$ and $E'\to T_{y_0}N$ have zero-dimensional fibers.
Their total spaces are therefore both naturally identified with $T_{y_0}N$; viewed over $M$, the total space of $E'$ is the vector space $T_{y_0}N$.
Thus the Fourier transform in (2) of Definition \ref{Pushforward}, which is taken along the zero-dimensional fibers, is the identity.
The Thom map in (1) is not the identity: it is $\maT_{T_{y_0}N}$ from the point to $T_{y_0}N$.
Also we have $M \times M \times N = N$ and the Morita equivalence in (4) is again the identity.

We turn to the map $(e _0)^{-1}$ in (3).
When $M$ is a single point $\maD$ is just the space
$$ \maD = T _{y_0} N \sqcup (N \times (0, 1]) \rightrightarrows T _{y _0} N \sqcup (N \times (0, 1]).$$
From the connection $\nabla$, one can define geodesics as curves with parallel tangent vector field.
Since if $t \mapsto \lambda (t)$ is a geodesic $t \mapsto \lambda (r t)$ is also a geodesic,
the exponential map $\exp _{y _0} : U \to N$ is well defined for some star-shaped open neighborhood $U$ of $0 \in T_{y _0} N$,
and is a diffeomorphism onto its image $V \subset N$.
Define a map $\maE : C^\infty _c (U) \to C^\infty _c (\maD)$ by
\begin{align*}
(\maE \phi ) (v) := \phi (v) & \quad \text{ if } v \in T _{y_0} N \\
(\maE \phi ) (y, r) := \phi (v') & \quad \text{ if } (y, r) \in N \times (0, 1], y = \exp _{y_0} ( r v' ), v' \in U \\
(\maE \phi ) (y, r) := 0 & \quad \text { if } y \not \in \exp _{y _0} r U.
\end{align*}
By uniqueness of geodesics, for any $v \in U, g \in G$, whenever $v \cdot g \in U$,
$$\exp _{y_0} (v \cdot g) = (\exp _{y_0} v) \cdot g .$$
It follows that $\maE$ extends to a homomorphism of convolution algebras
$$ \maE : \maC^\infty _c (U \rtimes G) \to \maC^\infty _c (\maD \rtimes G).$$
For any $\phi \in \maC^\infty _c (U \rtimes G)$, $e _0 \circ \maE (\phi)$ is $\phi$,
extended by zero to $\maC^\infty _c (T_{y_0}N \rtimes G)$.
By the Thom isomorphism, such extension induces a quasi-isomorphism in groupoid homology, hence periodic cyclic homology,
and therefore $\maE = (e _0)^{-1}$.
Evaluating $\maE \phi $ at $r = 1$ for any $\phi \in C^\infty _c (U)$, we obtain
$$ e _1 \circ (e _0 )^{-1} = \varepsilon _{N \rtimes G} \circ (\exp _{y_0} ^{-1})^*, $$
where $\varepsilon _{N \rtimes G} : \maC^\infty _c (V \rtimes G) \to \maC^\infty_c (N \rtimes G)$ is the extension by $0$ map.

To conclude, when $M$ is a point the pushforward map is given by

\begin{equation}
\label{PointPush}
f _! = \varepsilon _{N \rtimes G} \circ (\exp _{y_0} ^{-1})^* \circ \maT _{T _{y_0} N} = e _1 \circ (e _0 )^{-1} \circ \maT _{T _{y_0} N}.
\end{equation}

\section{Essential equivalence of \'etale Lie groupoids}
In this section we recall some properties that are invariant under essential
equivalence of (\'etale) Lie groupoids, namely the isotropy subgroup at a point,
the inertia groupoid and the groupoid homology.

\begin{dfn}
\label{EssEquivDfn}
Let $\varphi : \maK \to \maG$ be a homomorphism of Lie groupoids.
We say that $\varphi$ is an essential equivalence if:
\begin{enumerate}
\item
The map
$$ s \circ \pi _{\maG ^{(1)}} : \maK ^{(0)} \times _{\maG ^{(0)}} \maG ^{(1)} \to \maG ^{(0)} $$
(where $\maK ^{(0)} \times _{\maG ^{(0)}} \maG ^{(1)}
:= \{ (x, \gamma) \in \maK ^{(0)} \times \maG ^{(1)} : \varphi ^{(0)} (x) = \bt (\gamma) \} $)
is a surjective submersion;
\item
The map
$$ (\varphi ^{(1)} , \bs , \bt ) : \maK ^{(1)}
\to \maG ^{(1)} \times _{\maG ^{(0)} \times \maG ^{(0)}} (\maK ^{(0)} \times \maK ^{(0)}) $$
is a diffeomorphism.
\end{enumerate}
\end{dfn}

If $\maK, \maG$ are \'etale Lie groupoids, we shall further assume that $s \circ \pi _{\maG ^{(1)}}$ is \'etale.
This implies $\dim \maK ^{(1)} = \dim \maK ^{(0)} = \dim \maG ^{(0)} = \dim \maG ^{(1)} $ and hence $\varphi$ is also \'etale.

\begin{dfn}
Given Lie groupoids $\maG, \maG'$, we say that $\maG$ is Morita equivalent to $\maG'$,
if there exist essential equivalences $\varphi : \maK \to \maG, \varphi' : \maK \to \maG'$.
Usually we denote $\maG \xleftarrow{\varphi} \maK \xrightarrow{\varphi'} \maG'$.
\end{dfn}

Thus essential equivalence is the one-sided form of Morita equivalence:
a Morita equivalence is represented by a common refinement through essential equivalences.
The following example illustrates this relation for transformation groupoids.
After passing to a finite-index normal subgroup, the quotient map gives the refinement realizing a Morita equivalence,
and the invariance results below identify the corresponding cyclic homology and localized index pairings.

\subsection{Example: Morita equivalence of transformation groupoids from lattices in Lie groups}

We construct a class of examples illustrating the Morita equivalence of transformation groupoids.
We focus on the proper action of a lattice $\Gamma$ (possibly non-cocompact) in a semisimple Lie group $G$
on its associated symmetric space. By passing to a finite-index normal subgroup,
we replace the intricate dynamics of the infinite discrete group $\Gamma$
with a Morita-equivalent groupoid representing a finite group acting on a finite-volume manifold.
This simplifies the evaluation of the localized index pairing.

\subsubsection{The general construction}

Let $G$ be a connected semisimple Lie group with finite center, and let $K \subset G$ be a maximal compact subgroup.
The quotient $X = G/K$ is a symmetric space of non-compact type. The group $G$ acts transitively on $X$,
and the stabilizer of any point $x \in X$ is a conjugate of the compact group $K$.

Let $\Gamma \subset G$ be a lattice, for example an arithmetic lattice. The action of $\Gamma$ on $X$ is proper,
but not necessarily free. Because $\Gamma$ is discrete and the stabilizers in $G$ are compact,
the stabilizer of any point $x \in X$ in $\Gamma$ is the intersection of a discrete set with a compact set,
which is always a finite group. Thus, all isotropy subgroups of the $\Gamma$-action on $X$ are finite.

To resolve the singularities of this action, we rely on a foundational result in geometric group theory.

\begin{lem}[Selberg's Lemma]
Every finitely generated linear group over a field of characteristic zero contains a torsion-free subgroup of finite index.
\end{lem}

Because $G$ is a linear Lie group, the lattice $\Gamma$ is a finitely generated linear group.
By Selberg's Lemma, there exists a torsion-free subgroup $\Gamma_1 \subset \Gamma$ of finite index.
To ensure this subgroup behaves well under quotients, we define $N$ to be the normal core of $\Gamma_1$ in $\Gamma$:
$$ N = \bigcap_{\gamma \in \Gamma/\Gamma_1} \gamma \Gamma_1 \gamma^{-1} .$$
Since the index $[\Gamma : \Gamma_1]$ is finite, this intersection consists of finitely many conjugates.
Consequently, $N$ is a normal, finite-index, torsion-free subgroup of $\Gamma$.

Because $N$ is torsion-free, its only finite subgroup is the trivial group.
Therefore, the finite isotropy subgroups of the $N$-action on $X$ must be trivial, meaning $N$ acts freely and properly on $X$.

The quotient space $M = X/N$ is a smooth, non-compact manifold of finite volume.
Because $N$ is a normal subgroup of $\Gamma$, the quotient group $F = \Gamma/N$ is a finite group that acts smoothly on $M = X/N$.
Geometrically, the quotient map $X \to X/N$ is an $N$-principal bundle. In the language of Lie groupoids,
this gives an essential equivalence from the transformation groupoid $X \rtimes \Gamma$ to the action groupoid $(X/N) \rtimes F$.
Indeed, the object map $X \to X/N$ is a surjective submersion,
and the arrow map identifies $X \rtimes \Gamma$ with the pullback of $(X/N) \rtimes F$ along $X \to X/N$.
Hence one obtains a canonical Morita equivalence
$$ X \rtimes \Gamma \simeq (X/N) \rtimes F .$$
This Morita equivalence allows us to study the index pairing localized at one point of the infinite group $\Gamma$ acting on $X$
by equivalently studying the finite group $F$ acting on the finite-volume manifold $M$.

\subsubsection{A precise example: $\mathrm{SL}(n, \mathbb{Z})$}

To illustrate this construction concretely, consider the special linear group $G = \mathrm{SL}(n, \mathbb{R})$ for $n \ge 2$,
with maximal compact subgroup $K = \mathrm{SO}(n)$.
The associated symmetric space is $X = \mathrm{SL}(n, \mathbb{R})/\mathrm{SO}(n)$,
which can be identified with the space of positive-definite symmetric $n \times n$ matrices of determinant one.

Let $\Gamma = \mathrm{SL}(n, \mathbb{Z})$. This is a classic example of a non-cocompact lattice in $G$.
The action of $\mathrm{SL}(n, \mathbb{Z})$ on $X$ has non-trivial finite stabilizers, corresponding to matrices of finite order.

Instead of relying abstractly on Selberg's Lemma,
we can construct the torsion-free normal subgroup $N$ explicitly using principal congruence subgroups.
Let $p \ge 3$ be a prime number, and define the level-$p$ principal congruence subgroup $\Gamma(p)$
as the kernel of the reduction modulo $p$ map:
$$ 1 \to \Gamma(p) \to \mathrm{SL}(n, \mathbb{Z}) \to \mathrm{SL}(n, \mathbb{Z}/p\mathbb{Z}) \to 1 .$$
By definition, $\Gamma(p)$ is a normal subgroup of finite index in $\mathrm{SL}(n, \mathbb{Z})$.

To see that $\Gamma(p)$ is torsion-free for $p \ge 3$,
suppose $A \in \Gamma(p)$ has finite order $k > 1$. Because $A \equiv I \pmod{p}$,
we can write $A = I + p^m B$ where $m \ge 1$ and $B$ is a matrix not divisible by $p$.
A standard binomial expansion of $A^k = I$ leads to a contradiction in $p$-adic valuation unless $p=2$ and $m=1$.
Hence, for $p \ge 3$, $\Gamma(p)$ contains no elements of finite order.

Setting $N = \Gamma(p)$, we see that $N$ acts freely and properly on the symmetric space $X$.
The quotient $M = X/\Gamma(p)$ is a smooth manifold of finite volume.
The finite quotient group $F = \mathrm{SL}(n, \mathbb{Z}/p\mathbb{Z})$ acts smoothly on $M$, yielding the explicit Morita equivalence:
$$ X \rtimes \mathrm{SL}(n, \mathbb{Z}) \simeq (X/\Gamma(p)) \rtimes \mathrm{SL}(n, \mathbb{Z}/p\mathbb{Z}) .$$
This packages the complicated arithmetic dynamics into a finite group acting on a smooth,
non-compact manifold, rendering localized index trace evaluations highly tractable.

\subsection{The inertia groupoid is invariant under essential equivalence}

In this subsection, we assume we are given an essential equivalence $\varphi : \maK \to \maG$.

For each $x \in \maK ^{(0)}$, denote the isotropy subgroup at $x$ by
$$ \maK ^x _x := \{ \kappa \in \maK ^{(1)} : \bt (\kappa) = \bs (\kappa) = x \} .$$

\begin{prop}
\label{IsotropySubgroup}
If $\varphi : \maK \to \maG$ is an essential equivalence,
$\varphi$ restricts to a group isomorphism from $\maK ^x _x$ to $\maG ^{\varphi (x)}_{\varphi (x)} $.
\end{prop}
\begin{proof}
We restrict the diffeomorphism $(\varphi , \bs _{\maK}, \bt _\maK)$ in (2) of Definition \ref{EssEquivDfn} to $\maK ^x_x$.
The image is
$$ \maG ^{\varphi (x)}_{\varphi (x)} \times _{\maG ^{(0)} \times \maG ^{(0)}} \{ (x, x) \}, $$
which is isomorphic to $\maG ^{\varphi (x)}_{\varphi (x)}$.
\end{proof}

We denote the loop space of $\maG$ and $\maK$ by $B ^{(0)}$ and $\bar B ^{(0)}$ respectively,
and the inertia groupoid by $\maZ$ (resp. $\bar \maZ$) (see Definition \ref{inertialgroupoid}).
Observe that $\varphi ^{(1)} : \maK ^{(1)} \to \maG ^{(1)}$ restricts to a map
from $\bar B ^{(0)} $ to $B ^{(0)}$.
Moreover, we have the pullback square
\begin{equation}
\label{LoopSquare}
\begin{CD}
\bar B ^{(0)} @>{\varphi ^{(1)}}>> B ^{(0)} \\
@V{\bs _\maK = \bt _\maK}VV @VV{\bs _\maG = \bt _\maG}V \\
\maK ^{(0)} @>{\varphi ^{(0)}}>> \maG ^{(0)}
\end{CD}
\end{equation}
Hence $\varphi $ induces a groupoid homomorphism
$\bar \varphi : \bar \maZ \to \maZ$ defined by
$$ \bar \varphi ^{(0)} := \varphi ^{(1)} \vert _{\bar B ^{(0)}},
\quad \bar \varphi ^{(1)} ( \bar x, \kappa ) := (\varphi ^{(1)} (\bar x), \varphi ^{(1)} (\kappa)) ,
\quad \forall \, (\bar x , \kappa ) \in \bar B ^{(0)} \times _\bt \maK ^{(1)}.$$

The following statement is formulated in \cite[Page 53]{Ruan;Book} (second paragraph).
We give a proof here.

\begin{prop}
The groupoid homomorphism $\bar \varphi : \bar \maZ \to \maZ$ is an essential equivalence.
\end{prop}

\begin{proof}
To verify axiom (1) of Definition \ref{EssEquivDfn}, suppose we are given any $y \in B ^{(0)}$.
Using the hypothesis that $\varphi$ is an essential equivalence,
there exists $x, \in \maK ^{(0)} , \gamma \in \maG ^{(1)}$ satisfying
$$ \varphi ^{(0)} (x) = \bt _\maG (\gamma), \quad \bs _\maG (\gamma) = \bt _\maG (y) \in \maG ^{(0)} .$$
By \eqref{LoopSquare} there exists $z \in \bar B ^{(0)} $ such that
$$ \bar \varphi ^{(0)} (z) = \varphi ^{(1)} (z) = \gamma y \gamma ^{-1}, \quad \bt _\maK (z) = x .$$
Clearly $z $ together with $(\gamma y \gamma ^{-1} , \gamma ) \in \maZ$ satisfies axiom (1), i.e.,
$$ \bt _\maZ (\gamma y \gamma ^{-1} , \gamma ) = \bar \varphi ^{(0)} (z),
\quad \bs _\maZ (\gamma y \gamma ^{-1} , \gamma ) = y .$$

To verify axiom (2), suppose we are given any $z, z' \in \bar B ^{(0)}, (y, \gamma) \in \maZ ^{(1)}$ such that
$$\bs _\maZ (y, \gamma) = \gamma ^{-1} y \gamma = \varphi ^{(1)} (z) , \quad  \bt _\maZ (y, \gamma) = y = \varphi ^{(1)} (z').$$
These relations imply
$$ \bs _\maG (\gamma) = \varphi ^{(0)} (\bs _\maK (z)), \quad \bt _\maG (\gamma) = \varphi ^{(0)} (\bt _\maK (z')) .$$
Because $\varphi$ is an essential equivalence, there exists a unique $\kappa \in \maK ^{(1)}$ that satisfies
$$ \varphi ^{(1)} (\kappa ) = \gamma , \quad \bs _\maK (\kappa ) = \bs _\maK (z), \quad \bt _\maK (\kappa) = \bt _\maK (z'). $$
Consider $\kappa ^{-1} z' \kappa \in \bar B^{(0)}$. We have
\begin{align*}
\varphi ^{(1)} (\kappa ^{-1} z' \kappa)
= \gamma ^{-1} \varphi ^{(1)} (z') \gamma
=& \gamma ^{-1} y \gamma = \varphi ^{(1)} (z) \\
\bs _\maK (\kappa ^{-1} z' \kappa)
= \bs _\maK (\kappa) =& \bs _\maK (z).
\end{align*}
It follows from \eqref{LoopSquare} that $\kappa ^{-1} z' \kappa = z$ and
$(z' , \kappa ) \in \bar \maZ $ is the unique element that satisfies
$$ \bar \varphi ^{(1)} (z', \kappa) = (y, \gamma),
\quad \bs _{\bar \maZ} (z', \kappa) = z, \quad \bt _{\bar \maZ} (z', \kappa ) = z' . \qedhere$$
\end{proof}

\subsection{Essential equivalence induces isomorphism of groupoid homology}
Let $\maK, \maG$ be \'etale Lie groupoids and $\varphi : \maK \to \maG$ be an essential equivalence.
Let $\maA'$ be a $\maG$-sheaf, then $\maA := \varphi ^* \maA'$ is a $\maK$-sheaf.

\begin{lem}
\label{HomMorita}
\cite[Theorem 4.4, Corollary 4.6]{Moerdijk;EtaleHomology}
There exists a canonical isomorphism
$$ H _\bullet (\maG , \maA') \cong H _\bullet (\maK, \maA) .$$
\end{lem}

The key to proving Lemma \ref{HomMorita} is the construction of the comma groupoid, which we recall here.

\begin{dfn}
Let $\varphi : \maK \to \maG$ be a groupoid homomorphism. Fix $y \in \maG ^{(0)}$.
The comma groupoid $y / \varphi $ has object set
$$ y / \varphi ^{(0)} := \{ (x, \gamma) \in \maK ^{(0)} \times \maG ^{(1)} : \bs (\gamma) = y, \bt (\gamma) = \varphi (x) \},$$
and the arrow set
$$ y / \varphi ^{(1)} := \{ (\kappa , \gamma ) \in \maK ^{(1)} \times \maG ^{(1)}
: \bs (\gamma) = y, \bt (\gamma ) = \bs (\varphi (\kappa)) \}.
$$
Equivalently, $y / \varphi ^{(0)}$ is the set of diagrams
$$ x \to \varphi (x) \xleftarrow{\gamma} y ;$$
While $y / \varphi ^{(1)}$ is the set of commutative diagrams
\begin{equation}
\begin{CD}
\bs (\kappa) @>{\varphi}>> \bt (\gamma) @<{\gamma}<< y \\
@V{\kappa}VV @V{\varphi (\kappa)}VV @| \\
\bt (\kappa) @>{\varphi}>> \bt (\gamma') @<{\gamma'}<< y \\
\end{CD}
\end{equation}
where $\gamma' :=  \varphi (\kappa) \gamma $.

The source and target maps are respectively defined by taking the top and bottom rows of the above diagram,
and multiplication defined by composing along a row.
\end{dfn}

More generally, if $S$ is a submanifold of $\maG ^{(0)}$ we define the comma groupoid by
\begin{align*}
S / \varphi ^{(0)} :=& \{ (x, \gamma) \in \maK ^{(0)} \times \bs^{-1} (S) : \bt (\gamma) = \varphi (x) \}, \\
S / \varphi ^{(1)} :=& \{ (\kappa , \gamma ) \in \maK ^{(1)} \times \bs^{-1} (S) : \bt (\gamma ) = \bs (\varphi (\kappa)) \}.
\end{align*}

We write down the comma groupoid explicitly when $\varphi$ is an essential equivalence.
By Equation \eqref{LoopSquare},
given any pair $ x \to \varphi (x) \xleftarrow{\gamma} y , x' \to \varphi (x') \xleftarrow{\gamma'} y \in y / \varphi ^{(0)}$,
there exists a unique $\kappa \in \maK ^{(1)}$ that fits into the commutative diagram
$$
\begin{CD}
\kappa @>{\varphi }>> \gamma' \gamma^{-1} \\
@V{(\bs , \bt)}VV @V{(\bs, \bt)}VV \\
(x, x') @>{(\varphi , \varphi )}>> (\bt (\gamma), \bt (\gamma'))
\end{CD}
$$
by (2) of Definition \ref{EssEquivDfn}.
It is straightforward to verify $(\kappa, \gamma) \in y / \varphi ^{(1)}$.
Hence we conclude with the following lemma.

\begin{lem}
Suppose $\varphi : \maK \to \maG$ is an essential equivalence, then $y / \varphi$ is the pair groupoid
$$ y / \varphi ^{(0)} \times y / \varphi ^{(0)}.$$
\end{lem}


In order to prove Proposition \ref{EmbeddingMorita} below,
we need to recall how the isomorphism in Lemma \ref{HomMorita} is constructed, following \cite[Section 4.2]{Moerdijk;EtaleHomology}.

To begin with we define the maps
\begin{align*}
\alpha _n : \maK ^{(n)} \times _{\maG ^{(0)}} \maG ^{(1)} \to \maK ^{(0)}
& \quad \alpha_n (\kappa_1, \cdots , \kappa_n, \gamma ) := \bt (\kappa _1) \\
\beta _n : \maK ^{(n)} \times _{\maG ^{(0)}} \maG ^{(1)} \to \maG ^{(0)}
& \quad \beta _n (\kappa_1, \cdots , \kappa_n, \gamma ) := \bs (\gamma)
\end{align*}
and consider the sheaves
$$ \maB _n := (\beta _n)_! \alpha _n^* \maA $$
over $\maG ^{(0)}$.
Moreover, $\maB _\bullet$ is endowed with the structure of a complex of sheaves \cite[p.35]{Moerdijk;EtaleHomology}.

The stalk of $\maB _n$ is
$$ (\maB _n )_y = \Gamma _c (\beta _n^{-1} (y), \alpha _n ^* \maA) = B _n (y / \varphi , \pi _y^* \maA ) ,$$
where $\pi _y : (y / \varphi )^{(0)} = \maK ^{(0)} \times _{\maG ^{(0)}} \bs ^{-1} (y) \to \maK ^{(0)}$
is the projection onto the $\maK ^{(0)}$ component.

When $\maA = \varphi ^* \maA'$,
the stalk of $\pi_y ^* \varphi ^* \maA'$ at $\kappa \to \varphi (\kappa) \xleftarrow{\gamma} y $ is $\maA' _{\bt (\gamma)}$,
which can be canonically identified with $\maA' _y$ by action of $\gamma$.
The boundary operator on $\maB _\bullet$, when restricted to each stalk,
is just the groupoid boundary operator on $B _n (y / \varphi , \pi _y^* \maA )$.
Hence the groupoid homology $H _\bullet (y / \varphi , \pi _y^* \maA )$ (which is the pair groupoid) is well known, namely:
\begin{align*}
H _n (y / \varphi , \pi _y^* \maA ) \cong & \{ 0 \}, \quad n \geq 1 \\
H _0 (y / \varphi , \pi _y^* \maA ) \cong & \maA' _y,
\end{align*}
where the last isomorphism is given by
$$ c \mapsto \sum _{\tilde x \in y / \varphi ^{(0)}} c (\tilde x) .$$
It follows that $\cdots \to \maB _1 \to \maB _0 \to \maA' \to 0$ is a resolution of sheaves,
and the groupoid homology
$H _\bullet (\maG , \maA') $ is computed by the total homology of the double complex $H _\bullet (\maG , \maB _\bullet) $.

On the other hand, there is another way to compute $H _\bullet (\maG , \maB _\bullet) $.
Observe that the complex computing $H _\bullet (\maG , \maB _\bullet) $ is
\begin{equation}
\label{HomBi}
\maC _{p, q} := \Gamma _c (\maG ^{(p)}, \bt^* _p (\beta _q )_! \alpha _q ^* \maA )
= \Gamma _c (\maK ^{(q)} , (\nu _{p, q} )_! \nu _{p, q} ^* \bt _q^* \maA ),
\end{equation}
where $\nu _{p, q} : \maK ^{(q)} \times _{\maG ^{(0)}} \maG ^{(p+1)} \to \maK ^{(q)}$
is the projection onto the $\maK ^{(q)}$ component.

Summation along fiber gives a map from $\Gamma _c (\maK ^{(q)} , (\nu _{0, q} )_! \nu _{0, q} ^* \bt _q^* \maA )$
to the bar chain complex $B _q (\maK ^{(q)} , \bt _q ^* \maA ) = \Gamma _c (\maK ^{(q)} , \bt _q^* \maA) $.
Indeed one has exact sequences
$$ \cdots \to \maC _{1, q} \to \maC _{0, q} \to B _q (\maK ^{(q)} , \bt _q ^* \maA ) \to 0,$$
which implies
$$H _\bullet (\maG , \maB _\bullet) \cong H _\bullet (\maK , \maA),$$
as asserted in Lemma \ref{HomMorita}.

\subsection{Cyclic homology is invariant under essential equivalence}
Let $\varphi : \maK \to \maG$ be an essential equivalence.
We construct the isomorphism
$$ H P _\bullet (\maC ^\infty _c (\maK)) \cong H P _\bullet (\maC ^\infty _c (\maG)), $$
by applying Lemma \ref{HomMorita} to groupoid homology of the inertia groupoids,
for the essential equivalence $\bar \varphi : \bar \maZ ^{(0)} \to \maZ ^{(0)}$.

Observe furthermore that $\bar \varphi ^{(0)}$ maps $\bar B ^{(0)} _\ZZ $ to $B ^{(0)} _\ZZ$,
and $\bar B ^{(0)} _\infty $ to $B ^{(0)} _\infty$,
hence restricts to essential equivalence
$\bar \varphi : \bar \maZ _\ZZ \to \maZ _\ZZ$ and respectively $\bar \varphi : \bar \maZ _\infty \to \maZ _\infty$.
Moreover $\bar \varphi$ commutes with $\theta$.

For the hyperbolic case, suppose both $\bar \maN _\infty$ and $\maN _\infty$ have finite cohomological dimension
(indeed the groupoid cohomology of $\bar \maN _\infty$ with coefficient in $\CC$ is isomorphic to that of $\maN _\infty$)
we have by Lemma \ref{FiniteDim}
$$ H P _\bullet (\maC (\bar \maZ _\infty , \theta , C^\infty _{\maK ^{(0)}} )) \cong 0
\cong H P _\bullet (\maC (\maZ _\infty , \theta , C ^\infty _{\maG ^{(0)}})) .$$

For the elliptic case, we have
$$ H P _\bullet (\maC (\bar \maZ _\ZZ , \theta , \maA))
\xrightarrow{\cong} \prod _{k} H _{\bullet + 2 k} (\bar \maZ _\ZZ , \CC)
\xrightarrow{\cong} \prod _{k} H _{\bullet + 2 k} (\maZ _\ZZ , \CC)
\xleftarrow{\cong} H P _\bullet (\maC (\maZ _\ZZ , \theta , \maA)) ,$$
using Lemmas \ref{Elliptic} and \ref{HomMorita}.
If furthermore $\bar \maN _\infty$ has finite cohomological dimension,
then
$$ H P _\bullet (\maC ^\infty _c (\maK)) \cong  H P _\bullet (\maC (\bar \maZ _\ZZ , \theta , C^\infty _{\maK ^{(0)}})) ,$$
and similar relation holds for $ H P _\bullet (\maC ^\infty _c (\maG)) $.

\section{Index pairing}
With the preliminaries in the last two sections,
we define our index pairing in this section and prove its invariance under essential/Morita equivalence.

\subsection{A map from $H P _\bullet (\maC^\infty _c (\maG ^y _y))$ to $H P _\bullet (\maC^\infty _c (\maG))$}
In this section, let $\maG$ be a proper \'etale Lie groupoid.
We assume that $B ^{(0)}$ is oriented and the arrows in
$\maZ ^{(1)}$ are (germs of) orientation preserving diffeomorphisms.
Fix a point $y \in \maG ^{(0)}$ and let $\maG ^y _y $ be the isotropy subgroup of $y$.
Since $\maG$ is proper, we always have:
\begin{enumerate}
\item
There is a fixed $\maG$-invariant connection $\nabla$ on $T \maG ^{(0)}$;
\item
There is a fixed star-shaped $G$-invariant neighborhood $U \subset T_y \maG ^{(0)}$ of $0$,
such that $\exp _y : U \to \maG ^{(0)}$ is well defined and is a diffeomorphism to its image $V \subset \maG ^{(0)}$;
\item
The exponential map $\exp _y$ lifts to a groupoid isomorphism
$$ \exp _y : U \rtimes \maG ^y_y \to \maG _V \subset \maG .$$
\end{enumerate}
The upshot is that we can find some $V$ such that
$$ \maG _V \cong V \rtimes \maG ^y _y .$$

Let $f _y$ be the embedding of single point space $\{ y \}$ into $V$.
Then $f _y$ is equivariant,
or equivalently $f _y$ induces a groupoid homomorphism from $G$ to $V \rtimes \maG ^y_y $.
Hence we can define $(f_y)_!$ as in Definition \ref{Pushforward}.

Observe that $\maG _V $ is an open sub-groupoid of $\maG$.
Moreover, from Definition \ref{Star}, for any $f_1, f_2 \in \maC _c ^\infty (\maG)$ supported inside $\maG _V$,
$f_1 (\gamma \gamma _1 ^{-1}) f_2 (\gamma _1) = 0$ unless $\bs (\gamma _1), \bt (\gamma _1), \bs (\gamma), \bt (\gamma) \in V$.
It follows that extension by zero is an algebra homomorphism from $\maC^\infty _c (\maG _V)$ to $\maC^\infty _c (\maG)$,
which we denote by $\varepsilon _\maG $.
We consider the induced map
$$ \varepsilon _\maG \circ f_! : H P _\bullet (\maC^\infty _c (\maG ^y _y))
\to H P _{\bullet + \dim \maG ^{(0)}} (\maC^\infty _c (\maG)).$$

Note that for any other star-shaped, $\maG ^y_y$-invariant neighborhood $U' \subset U$ of $0$ and $V' := \exp _y U' \subset V$,
one has the commutative diagram
$$
\begin{CD}
H P _\bullet (\maC^\infty _c (\maG ^y _y)) @>{(f _y)_!}>>
H P _{\bullet + \dim \maG ^{(0)}} (\maC^\infty _c (\maG _{V'}))
@>{\varepsilon _\maG }>>H P _{\bullet + \dim \maG ^{(0)}} (\maC^\infty _c (\maG)) \\
@| @V{\varepsilon _{V', V}}VV @| \\
H P _\bullet (\maC^\infty _c (\maG ^y _y)) @>{(f _y)_!}>>
H P _{\bullet + \dim \maG ^{(0)}} (\maC^\infty _c (\maG _{V}))
@>{\varepsilon _\maG }>>H P _{\bullet + \dim \maG ^{(0)}} (\maC^\infty _c (\maG))
\end{CD}
$$
where $\varepsilon _{V', V}$ is the extension by $0$ map,
which induces an isomorphism because $U'$ and $U$ are contractible
(see Example \ref{ThomEx} for explicit computation on the level of homology).
Therefore $ \varepsilon _\maG \circ f_! $ does not depend on the choice of $U$ or $V$.
Moreover since the connection $\nabla$ is only used to construct $V$ via the exponential map,
by similar shrinking arguments the construction does not depend on $\nabla$.

\begin{rem}
The orientation assumption can be weakened, as we only need the pushforward for the local subgroupoid $\maG _V$.
Also as pointed out in Remark \ref{NonProper},
$\varepsilon _\maG \circ f_! $ is well defined assuming the existence of a connection $\nabla$ as above,
instead of assuming $\maG$ to be proper.
\end{rem}

The main theorem of this section is that such $\varepsilon _\maG \circ (f _y)_!$ is compatible with essential equivalence.
Hence let $\varphi : \maK \to \maG $ be an essential equivalence that is orientation preserving on the loop spaces.
We fix $x \in \maK ^{(0)}$ and let $y := \varphi ^{(0)} (x) \in \maG ^{(0)}$.
We construct the connection $\bar \nabla$ and $\bar U \subset T_x \maK ^{(0)}$ as above.
We may furthermore assume that
\begin{enumerate}
\item
The connection $\bar \nabla$ on $T \maK ^{(0)}$ is the pullback of $\nabla$, via the differential $d \varphi$, which implies
$ \exp _y \circ (d \varphi) = \varphi ^{(0)} \circ \exp _x : \bar U \to V $;
\item
The differential $d \varphi ^{(0)}$ is a bijective map from $\bar U$ to $U$.
\item
The local diffeomorphism $\varphi ^{(0)}$ restricts to an equivariant diffeomorphism
$ \bar V := \exp _x (\bar U) \to V := \exp _y (U)$,
and hence $\varphi : \maK _{\bar V} \to \maG _V $ is a groupoid isomorphism,
and that
\begin{equation}
\label{CD0}
\begin{CD}
\bar U \rtimes \maK ^x_x @>>> \maK _{\bar V} @>>> \maK \\
@V{d \varphi}VV @V{\varphi}VV @V{\varphi}VV \\
U \rtimes \maG ^y_y @>>> \maG _V @>>> \maG
\end{CD}
\end{equation}
commutes.
\end{enumerate}

Here we make a further simplifying observation.
Since $\maK$ (and hence $\maG$) is proper, the hyperbolic part does not contribute to the periodic cyclic homology, i.e.,
$$ H P _\bullet (\maC (\bar \maZ _\infty , \theta , C^\infty _{\maK ^{(0)}} )) \cong 0
\cong H P _\bullet (\maC (\maZ _\infty , \theta , C ^\infty _{\maG ^{(0)}})) $$
(see Lemma \ref{FiniteDim}).
Thus we only need to consider the elliptic part of the periodic cyclic homology,
which can be computed by groupoid homology using Lemma \ref{Elliptic}.

\begin{prop}
We have the commutative diagram
\begin{equation}
\label{CD1}
\begin{CD}
H P _\bullet (\maC^\infty _c (\maK ^x _x))
@>{(f _x)_!}>> H P _{\bullet + \dim \maK ^{(0)}} (\maC^\infty _c (\bar V \rtimes \maK ^x _x)) \\
@V{(\varphi ^{-1})^*}VV @V{(\varphi ^{-1})^*}VV \\
H P _\bullet (\maC^\infty _c (\maG ^y _y))
@>{(f _y)_!}>> H P _{\bullet + \dim \maG ^{(0)}} (\maC^\infty _c (V \rtimes \maG ^y _y))
\end{CD}
\end{equation}
\end{prop}
\begin{proof}
Equation \eqref{CD1} follows from the functoriality Lemma \ref{PushLem1},
and the observation that $\varphi _! = (\varphi ^{-1})^*$ when $\varphi$ is a groupoid isomorphism.
\end{proof}

Alternatively, one can combine \eqref{CD0} with Thom isomorphism to get the commutative diagram
$$
\begin{CD}
\maK ^x_x @>{\maT _{T_0 \bar U}}>> \bar U \rtimes \maK ^x_x @>{\exp _x}>> \maK _{\bar V} \\
@V{\varphi}VV @V{d \varphi}VV @V{\varphi}VV \\
\maK ^y_y @>{\maT _{T_0 U}}>> U \rtimes \maG ^y_y @>{\exp _y}>> \maG _V
\end{CD}
$$
and \eqref{CD1} follows because $(f _y)_! = (\exp ^{-1} _y)^* \circ \maT _{T _y V}$ by Equation \eqref{PointPush}
(there is no need to extend by zero because $\exp _y : U \to V$ is surjective).

\begin{rem}
In the rest of this section, we only need the commutative diagram
$$
\begin{CD}
\maK ^x_x @>>> \maK _{\bar V} @>>> \maK \\
@V{\varphi}V{\cong}V @V{\varphi}V{\cong}V @V{\varphi}VV \\
\maG ^y_y @>>> \maG _V @>>> \maG
\end{CD}
$$
to hold for some contractible neighborhoods $V, \bar V$.
\end{rem}

Denote by $B ^{(0)} _V $ the loop space of $\maG _V$, and $\maZ _V $ the inertia groupoid of $\maG _V$.
Note that by definition
\begin{align*}
B ^{(0)} _V := & \{ \gamma \in \maG ^{(1)} : \bs (\gamma ) = \bt (\gamma ) \in V \} \\
\maZ ^{(1)} _V := & B ^{(0)} _V \rtimes \maG _V = \bt ^{-1} _{\maZ} (B ^{(0)} _V) \cap \bs ^{-1} _\maZ (B ^{(0)} _V).
\end{align*}
Similarly we denote by $B ^{(0)} _{\bar V} $ and $\bar \maZ _{\bar V}$,
the loop space and, respectively, the inertia groupoid of $\maK _{\bar V}$.

\begin{prop}
We have the commutative diagram
\begin{equation}
\label{CD2}
\begin{CD}
H P _\bullet (\maC^\infty _c (\maK _{\bar V} ))
@>>> H P _\bullet (\maC ((\bar \maZ _{\bar V})_\ZZ , \theta , \maA))
@>>> \prod _k H _{\bullet + 2k} ((\bar \maZ _{\bar V}) _\ZZ , \CC) \\
@V{(\varphi ^{-1})^*}VV @V{(\bar \varphi ^{-1})^*}VV @V{\bar \varphi _*}VV \\
H P _\bullet (\maC^\infty _c (\maG _{V} ))
@>>> H P _\bullet (\maC ((\maZ _V)_\ZZ , \theta , \maA))
@>>>\prod _k H _{\bullet +2k} ((\maZ _{V}) _\ZZ , \CC)
\end{CD}
\end{equation}
where the rows of the left square are the isomorphism in Lemma \ref{Reduction},
followed by projection onto the elliptic part,
and the rows of the right square is the isomorphism in Lemma \ref{Elliptic}.
\end{prop}

\begin{proof}
Since $\varphi : \maK _{\bar V} \to \maG _V$ is a groupoid isomorphism,
pulling back by $\varphi ^{-1}$ is clearly compatible with the canonical isomorphism of Lemma \ref{Reduction}.
Hence we get the commutative diagram:
$$
\begin{CD}
H P _\bullet (\maC^\infty _c (\maK _{\bar V} )) @>>> H P _\bullet (\maC (\bar \maZ _{\bar V} , \theta , \maA)) \\
@V{ (\varphi ^{-1})^*}VV @V{ (\bar \varphi ^{-1})^*}VV \\
H P _\bullet (\maC^\infty _c (\maG _{V} )) @>>> H P _\bullet (\maC (\maZ _V , \theta , \maA))
\end{CD}
$$
moreover $\bar \varphi$ preserves the elliptic/hyperbolic parts, therefore commutes with projection
$$
\begin{CD}
H P _\bullet (\maC ((\bar \maZ _{\bar V})_\ZZ , \theta , \maA))
\oplus  H P _\bullet (\maC ((\bar \maZ _{\bar V})_\infty , \theta , \maA))
@>>> H P _\bullet (\maC ((\bar \maZ _{\bar V})_\ZZ , \theta , \maA)) \\
@V{(\bar \varphi ^{-1})^*}VV @V{(\bar \varphi ^{-1})^*}VV \\
H P _\bullet (\maC ((\maZ _{V})_\ZZ , \theta , \maA))
\oplus  H P _\bullet (\maC ((\maZ _{V})_\infty , \theta , \maA))
@>>> H P _\bullet (\maC ((\maZ _V)_\ZZ , \theta , \maA))
\end{CD}
$$
and the left square follows.

Similarly the isomorphism in Lemma \ref{Elliptic} is canonical, and we get the commutative diagram
$$
\begin{CD}
H P _\bullet (\maC ((\bar \maZ _{\bar V})_\ZZ , \theta , \maA))
@>>> \prod _k H _{\bullet + 2k} ((\bar \maZ _{\bar V}) _\ZZ , \CC) \\
@V{(\bar \varphi ^{-1})^*}VV @V{(\bar \varphi ^{-1})^*}VV \\
H P _\bullet (\maC ((\maZ _V)_\ZZ , \theta , \maA))
@>>>\prod _k H _{\bullet +2k} ((\maZ _{V}) _\ZZ , \CC)
\end{CD}
$$
and since $\bar \varphi _*$ is by definition summation over the fibers (see Lemma \ref{HomFunct}),
when $\varphi $ is bijective, we have $\bar \varphi _* = (\bar \varphi ^{-1})^* $.
\end{proof}

Since $\maZ _V$ is an open subgroupoid of $\maZ $,
denote the embedding from $\maZ _{V} $ to $\maZ$, by $\iota$.
Then the induced map $\iota _* : B _n (\maZ _V , \Omega ^q ) \to B _n (\maZ , \Omega ^q )$ is just the extension by zero map.

\begin{prop}
\label{EmbeddingMorita}
We have the commutative diagram
\begin{equation}
\label{CD3}
\begin{CD}
H _{\bullet } ((\bar \maZ _{\bar V}) _\ZZ , \CC) @>{\iota _*}>> H _{\bullet } (\bar \maZ _\ZZ , \CC) \\
@V{\bar \varphi _*}VV @V{\cong}VV \\
H _{\bullet } ((\maZ _{V}) _\ZZ , \CC) @>{\iota _*}>> H _{\bullet } (\maZ _\ZZ , \CC)
\end{CD}
\end{equation}
where the map at the right column is the isomorphism induced by Morita equivalence (Lemma \ref{HomMorita}).
\end{prop}
\begin{proof}
Recall that by assumption $\varphi : (\bar \maZ _{\bar V}) _\ZZ \to (\maZ _{V}) _\ZZ$
is an isomorphism, hence also an essential equivalence.
Hence taking the c-soft resolution by sheaf of germs of differential forms,
we can repeat the arguments in the proof of Lemma \ref{HomMorita}, in particular constructing the bi-complexes
$ \maC _{p, q}$ in Equation \eqref{HomBi}, for the essential equivalence $\bar \maZ _\ZZ \xrightarrow{\bar \varphi} \maZ _\ZZ$,
and also $\tilde \maC _{p, q}$, with $(\bar \maZ _{\bar V}) _\ZZ \xrightarrow{\bar \varphi} (\maZ _{V}) _\ZZ$ in place of
$\maK \xrightarrow{\varphi} \maG$,
and coefficient sheaves $\Omega $ and $\iota ^* \Omega$ respectively.

Clearly $\iota _* : \tilde \maC _{p, q} \to \maC _{p, q} $ is a morphism of bicomplex
(which commutes with the differential maps).
Hence we have the commutative diagram
$$
\begin{CD}
H _\bullet ((\bar \maZ _{\bar V}) _\ZZ , \bar \varphi ^* \iota ^* \Omega)
@<<< H _\bullet (\oplus _{p+q = \bullet} \tilde \maC _{p, q} , \partial + \delta)
@>>> H _\bullet ( (\maZ _V ) _\ZZ , \iota ^* \Omega) \\
@V{\iota _*}VV @V{\iota _*}VV @V{\iota _*}VV \\
H _\bullet (\bar \maZ _\ZZ , \bar \varphi ^* \Omega)
@<<< H _\bullet (\oplus _{p+q = \bullet} \maC _{p, q} , \partial + \delta)
@>>> H _\bullet (\maZ _\ZZ , \Omega)
\end{CD}
$$
On the other hand, since $\bar \varphi : \bar \maZ _{\bar V} \to \maZ _{V}$ is an isomorphism,
it follows from \cite[Remark 4.5.3]{Moerdijk;EtaleHomology} that the isomorphism of groupoid homology induced by $\tilde \maC _{p, q}$
in the top row coincides with \cite[Section 3.7.3]{Moerdijk;EtaleHomology}
$$\bar \varphi _* : H _\bullet ( (\bar \maZ _{\bar V}) _\ZZ , \bar \varphi ^* \iota ^* \Omega)
\to H _\bullet ( (\maZ _V ) _\ZZ , \bar \varphi _! \bar \varphi ^* \iota ^* \Omega)
= H _\bullet ( (\maZ _V ) _\ZZ  , \iota ^* \Omega) . $$
The proposition follows because $\iota ^* \Omega$ is just the sheaf of germs of differential forms on $(B ^{(0)} _V ) _\ZZ $
(and respectively $(\bar B ^{(0)} _{\bar V} ) _\ZZ $), which are both resolutions of $\CC$.
\end{proof}

Define the isomorphism
$\maM _\varphi : H P _\bullet (\maC^\infty _c (\maK)) \to H P _\bullet (\maC^\infty _c (\maG)) $ by the commutative diagram
\begin{equation}
\label{CD4}
\begin{CD}
\prod _{k} H_ {\bullet + 2 k} (\bar \maZ _\ZZ , \CC) @<<< H P _\bullet (\maC^\infty _c (\maK)) \\
@VVV @V{\maM _\varphi}VV \\
\prod _{k} H_ {\bullet + 2 k} (\maZ _\ZZ , \CC) @<<< H P _\bullet (\maC^\infty _c (\maG))
\end{CD}
\end{equation}
where the map at the left column is the isomorphism induced by Morita equivalence (Lemma \ref{HomMorita}),
as in \eqref{CD3}.
Hence we consider the composition
\begin{align}
\label{TopComp}
H P _\bullet & (\maC^\infty _c (\maG _V)) \to H P _\bullet (\maC ((\maZ _{V})_\ZZ , \theta , \maA))
\to \prod _k H _{\bullet + 2 k} ( (\maZ _V ) _\ZZ , \CC) \\ \nonumber
\xrightarrow{\iota _*} & \prod _{k} H_ {\bullet + 2 k} (\maZ _\ZZ , \CC) \to H P _\bullet (\maC ((\maZ )_\ZZ , \theta , \maA))
\to H P _\bullet (\maC^\infty _c (\maG)) .
\end{align}
Note that the first line is just the bottom row of \eqref{CD2}.

\begin{prop}
(Compare \cite[Theorem 4.4]{Crainic;EtaleGpoid}) We have the commutative diagram
\begin{equation}
\label{CD5}
\begin{CD}
H P _\bullet (\maC^\infty _c (\maG _V)) @>{\varepsilon _\maG}>> H P _\bullet (\maC^\infty _c (\maG)) \\
@VVV @VVV \\
\prod _k  H _{\bullet + 2 k} ( (\maZ _V ) _\ZZ , \CC)
@>{\iota _*}>> \prod _{k} H_ {\bullet + 2 k} (\maZ _\ZZ , \CC)
\end{CD}
\end{equation}
where the left column (respectively, right column) is the first line of \eqref{TopComp} (resp. the inverse of the second line),
and similar for $\maK$.
\end{prop}
\begin{proof}
It is straightforward to verify that the construction leading to Lemma \ref{Elliptic} is compatible with extension by zero.
Here we just recall the construction steps and sketch the compatibility arguments.
Again we shall write $\maA := C^\infty _{\maG ^{(0)}}$ as in Example \ref{CrossProd},
and we denote $\maA _V := C^\infty _{\maG ^{(0)} _V} = C^\infty _V$.
Note that $\maA _V$ is just the restriction of $\maA$ to $V$.

To begin with, define $(\maA \rtimes \maG) ^\natural $ as in \eqref{AlgCross}.
One can repeat the same construction for $\maG _V$, and define $B ^{(n)}_V$ and $(\maA _V \rtimes \maG _V) ^\natural$.

Recall that in Lemma \ref{Reduction} the quasi-isomorphism from
$\maC^{\infty}_c (\maG ) ^\natural$ to $(\maA \rtimes \maG) ^\natural$,
is induced by restricting $c$ to its sheaf over $B ^{(n)}$
(and respectively from $\maC^{\infty}_c (\maG _V )^\natural$ to $(\maA _V \rtimes \maG _V) ^\natural$).
Observe that $\maG _V $ is an open subset of $\maG$ and $B ^{(n)}_V = \maG ^{(n+1)} _V \cap B ^{(n)}$.
It follows that restriction commutes with extension by $0$.

Next, one identifies $(\maA \rtimes \maG) ^\natural$ with $\maC (\maZ , \theta , \maA)$ (see Definition \ref{TwistedLoop}),
by the bijection \eqref{Twist}.
Restricting \eqref{Twist} to $\gamma _i \in \maG _V$ and $z_j \in \maZ _V$ identifies
$(\maA _V \rtimes \maG _V) ^\natural$ with $\maC (\maZ _V , \theta , \maA _V)$.
Moreover both the elliptic and hyperbolic parts are open subsets,
therefore the extension by $0$ map restricts to each of them,
and we have the commutative diagram
$$
\begin{CD}
H P _\bullet (\maC ((\maZ _{V})_\ZZ , \theta , \maA))
\oplus  H P _\bullet (\maC ((\maZ _{V})_\infty , \theta , \maA))
@>>> H P _\bullet (\maC ((\maZ _{V})_\ZZ , \theta , \maA)) \\
@VVV @VVV \\
H P _\bullet (\maC (\maZ_\ZZ , \theta , \maA))
\oplus  H P _\bullet (\maC (\maZ _\infty , \theta , \maA))
@>{\cong}>> H P _\bullet (\maC (\maZ _\ZZ , \theta , \maA))
\end{CD}
$$
where the columns are extension by $0$ and the rows are projection.

Lastly, one applies the Eilenberg-Zilber theorem (which is canonical),
and the Hochschild-Kostant-Rosenberg quasi-isomorphism, which is local (see \cite[Lemma 3.5]{Crainic;EtaleGpoid}),
to compute $H P _\bullet (\maC (\maZ _\ZZ , \theta , \maA))$
in terms of groupoid homology
$H _\bullet (\maN _\ZZ , \CC) \cong H _\bullet (\maZ _\ZZ , \CC)$, as recalled in Lemma \ref{Elliptic}.
Therefore this quasi-isomorphism is again compatible with extension by $0$,
which is just the pushforward $\iota _*$ because $\iota : \maZ _V \to \maZ$ is injective.
\end{proof}

Combining the commutative diagrams \eqref{CD1}, \eqref{CD2}, \eqref{CD3}, \eqref{CD4},
and then computing the top and bottom rows with \eqref{CD5}, we conclude with the main theorem in this section.

\begin{thm}
\label{CD6}
We have the commutative diagram
$$
\begin{CD}
H P _\bullet (\maC^\infty _c (\maK ^x _x))
@>{\varepsilon _\maK \circ (f _x)_!}>> H P _{\bullet + \dim \maK ^{(0)}} (\maC^\infty _c (\maK)) \\
@V{(\varphi ^{-1})^*}VV @VVV \\
H P _\bullet (\maC^\infty _c (\maG ^y _y))
@>{\varepsilon _\maG \circ (f _y)_!}>> H P _{\bullet + \dim \maG ^{(0)}} (\maC^\infty _c (\maG))
\end{CD}
$$
\end{thm}

\subsection{Definition of index pairing}
Having constructed the map $\varepsilon \circ f_!$ from
$H P _\bullet (\maC^\infty _c (\maG ^y _y))$ to
$H P _{\bullet + \dim \maG ^{(0)}} (\maC^\infty _c (\maG))$ for each $y \in \maG ^{(0)}$,
we can readily define a pairing with cyclic cohomology,
$H P ^{\bullet + \dim \maG ^{(0)}}  (\maC^\infty _c (\maG))$.
Chern character of elliptic operators or $K$-homology class defines classes in cyclic cohomology.

\begin{dfn}
Given any $[c] \in H P _\bullet (\maC^\infty _c (\maG ^y _y))$,
and any $\Psi \in H P ^{\bullet + \dim \maG ^{(0)}} (\maC^\infty _c (\maG))$,
define the pairing
$$ \langle \varepsilon \circ (f _y)_! [c] , \Psi \rangle \in \CC,$$
where $\langle \cdot , \cdot \rangle $ denotes the pairing between periodic cyclic homology and cohomology of $\maC^\infty _c (\maG)$.
\end{dfn}

There is a canonical isomorphism $\maM _\varphi ^* : H P ^\bullet (\maC^\infty _c (\maK)) \to H P ^\bullet (\maC^\infty _c (\maG)) $
such that
$$ \langle \maM _\varphi ([C]) , \maM _\varphi ^* (\Psi) \rangle = \langle [C], \Psi \rangle $$
for any $[C] \in H P _\bullet (\maC^\infty _c (\maK)), \Psi \in H P ^\bullet (\maC^\infty _c (\maK)) $
\cite[Theorem 4.12]{Crainic;EtaleGpoid}.
Indeed, recall that periodic cyclic cohomology is computed by groupoid cohomology and one has Poincare duality in the \'etale case
\cite{Moerdijk;EtaleHomology, Behrend;CohStack}.
The following corollary follows immediately from Theorem \ref{CD6}.

\begin{cor}
\label{CD7}
Let $\varphi : \maK \to \maG$ be an essential equivalence. Fix $x \in \maK ^{(0)}$ and let $y := \varphi ^{(0)} (x)$.
For any $[c] \in H P _\bullet (\maC^\infty _c (\maG ^y _y))$, $\Psi \in H P ^{\bullet + \dim \maG ^{(0)}} (\maC^\infty _c (\maG))$,
$$ \langle \varepsilon \circ (f _y)_! ([c]) ,  \Psi \rangle
= \langle \varepsilon \circ (f _x)_! (\varphi ^* [c]) , (\maM _\varphi ^* )^{-1} \Psi \rangle .$$
\end{cor}

Corollary \ref{CD7} generalizes to Morita equivalence.
Suppose we are given a Morita equivalence $\maG \xleftarrow{\varphi} \maK \xrightarrow{\varphi '} \maG'$,
where $\varphi , \varphi'$ are essential equivalences (and hence restrict to isomorphisms of the respective isotropy subgroups).
\begin{cor}
Let $x \in \maK , y \in \maG , y' \in \maG'$ be such that $y = \varphi ^{(0)} (x), y' = \varphi ^{\prime (0)} (x)$.
Then for any $[c] \in H P _\bullet (\maC^\infty _c (\maG ^y _y))$
and $\Psi \in H P ^{\bullet + \dim \maG ^{(0)}} (\maC^\infty _c (\maG))$,
we have
$$ \big\langle \varepsilon \circ (f _y)_! ([c]) ,  \Psi \big\rangle
= \big\langle \varepsilon \circ (f _{y'})_! \big(((\varphi')^{-1})^* \varphi ^* [c] \big) ,
\maM ^* _{\varphi '} \circ (\maM ^* _\varphi )^{-1} \Psi \big\rangle ,$$
where $\varphi '$ in the first pairing on the right is the restriction $\varphi' : \maK ^x _x \to \maG ^{y'} _{y'}$,
which is a group isomorphism.
\end{cor}

\subsection{Conclusion and open questions}
We have constructed a pairing
$$ ([c], \Psi) \mapsto \langle \varepsilon \circ (f _y)_! ([c]) ,  \Psi \rangle \in \CC,$$
for any $[c] \in H P _\bullet (\maC ^\infty _c (\maG ^y_y)), \Psi \in H P ^{\bullet + \dim \maG ^{(0)}} (\maC ^\infty _c (\maG))$,
which in the special case of action groupoid with discrete fixed points,
computes the equivariant index pairing.
It seems that one can construct a similar but more general type of index pairing when the
(local) fixed point sets have positive dimension,
by considering appropriate stratification of local group actions.

\section{Relation to Atiyah's Distribution Index}
In this section, we explain how our index pairing recovers components of the higher index
for a finite symmetry group, where the higher index is identified with the classical equivariant index.
Let $M$ be a closed, oriented, even-dimensional spin manifold of dimension $n$,
and let $G$ be a finite group acting properly and isometrically on $M$, preserving the orientation and the spin structure.
Let $D ^+$ be the positive chiral part of the associated $G$-invariant spin Dirac operator $D$ on $M$.

For simplicity, we first assume that the set of points with non-trivial isotropy is the finite set
$\{x_1, \ldots, x_N\}$, and moreover,
there is a fixed subgroup $H \subset G$ such that $H$ is the isotropy subgroup of all $x_1, \ldots, x_N$.

Since $G$ is finite, $C^*_r(G) \cong C^*(G) \cong \mathbb{C}G$ and
$K _0(C^*(G)) \cong R(G)$.
The higher index of $D$ is therefore
$$\Ind_G(D)=[\ker D ^+]-[\operatorname{coker}D ^+]
\in K _0(C^*(G)) \cong R(G),$$
and under this identification it is the classical equivariant index.
Its character extends to a trace on $\mathbb{C}G$: for any $f \in \mathbb{C}G$,
$$[\Ind_G(D)](f) := \operatorname{Tr}(f|_{\ker D ^+})
- \operatorname{Tr}(f|_{\operatorname{coker}D ^+}) \in \mathbb{C}.$$
In particular, evaluating at the delta function $\delta_g$ gives the $g$-component,
or equivariant Lefschetz number,
$$\ind _g(D):=[\Ind_G(D)](\delta_g) \in \mathbb{C}.$$
The class $\Ind_G(D)$ is the pushforward of the $G$-equivariant $K$-homology class
of $D$ along the constant map $p: M \to pt$.

By Ponge-Wang \cite[Theorem 11.4]{PW18}, the class of $\ch (D)$ is represented by the current
\begin{equation}
\label{PWCh}
\Upsilon^g (R)\vert_{M_a^g}
= (-i)^{\frac{n}{2}}(2\pi)^{-\frac{a}{2}} \widehat{A}(R^{TM_a^g}) \wedge \mathcal{V}_g(R^{\mathcal{N}_a^g}),
\end{equation}
with $\mathcal{V}_g(R^{\mathcal{N}_a^g}) = \det_\CC ^{-\frac{1}{2}}\big(1 - g^{\mathcal{N}} e^{-R^{\mathcal{N}_a^g}} \big)$.
Since for $g \in H \setminus \{e\}$ the fixed-point set $M^g$ is zero dimensional by assumption,
\begin{equation}
\label{PWIsolated}
\Upsilon^g(R)(x_i,g)
=(-i)^{\frac{n}{2}}\det {}_\CC^{-\frac{1}{2}}(1-g^{\mathcal N}_{x_i}).
\end{equation}
The equivariant fixed-point formula (see \cite[Lemma 6.4]{BGVBook04}) therefore gives
\begin{equation}
\label{AB}
\ind_g(D)=\sum_{i=1}^N \Upsilon^g(R)(x_i,g)
\end{equation}
for any $g \in H \setminus \{e\}$, while $\ind_g(D)=0$ for $g \in G \setminus H$.

Let $f_i : \{ x_i \} \to M$ be the embedding corresponding to $x_i$,
and $(f_i)_! := (f _{x_i})_!$ be defined as in Definition \ref{Pushforward}.
Fix $l \in H \setminus \{e\}$ and let $\langle l \rangle _H $ be the conjugacy class of $l$ in $H$.
We consider elements in $H P _0 (\CC H)$ of the form
$$ [c] = \sum _{g \in \langle l \rangle _H } g $$
(because these elements span $H P _0 (\CC H)$).
Then $\varepsilon (f _i)_! [c] $, when identified with groupoid homology over $\maZ$ (as in Corollary \ref{HPSplitting}),
is represented by the class that is $1$ on $ (x_i , g) \in B ^{(0)}, g \in \langle l \rangle _H $, and $0$ otherwise.

\begin{thm}
We have
\begin{align*}
\sum _{i=1}^N \langle \varepsilon (f _i)_! [c] , \ch (D) \rangle
=& \sum _{g \in \langle l \rangle _H} \ind _g (D) 
\end{align*}
\end{thm}
\begin{proof}
Let $\langle l \rangle _H $ be the conjugacy class of $l$ in $H$, and
$[c] = \sum _{g \in \langle l \rangle _H } g \in H P _0 (\CC H)$ as above.
We compute $\langle \varepsilon (f _i)_! [c] , \ch (D) \rangle$ for an arbitrary fixed $f_i$.

As in Section 4, we fix an $H$-invariant star shaped neighborhood $V _i = \exp _{x_i} U_i \subset M$ of $x_i$.
By shrinking $V _i$,
we may assume $x_i$ is the only point with non-trivial subgroup in $V_i$,
and moreover for any $g \not \in H, g V_i \cap V _i = \emptyset$.
Then $(M \rtimes G) _{V_i} = V _i \rtimes H$, and the loop space of $V_i \rtimes H$ is
$$ B ^{(0)} _{V_i \rtimes H} = V_i \sqcup  (x_i \times (H \setminus \{ e \})), $$
where $H$ acts on $x_i \times (H \setminus \{ e \})$ by conjugating the $H \setminus \{ e \}$ component.
Recall furthermore that one has splitting
$$ x_i \times (H \setminus \{ e \})
= \sqcup _{\langle g \rangle _H \in \langle H \rangle \setminus \{ e \} } x_i \times \langle g \rangle .$$
Comparing with the action groupoid $\{ e \} \rtimes H = H$ and following the computation in Sections 2.5 and 2.6,
the pushforward between groupoid homology
$$(f_i)_! : H _\bullet (Z ^g , \CC) \cong H _\bullet ( \langle g \rangle _H \rtimes H , \CC)
\to H _\bullet ( (\{ x_i \} \times \langle g \rangle _H) \rtimes H , \CC) $$
is just the identity map for $g \neq e$.
Hence we conclude that in groupoid homology,
\begin{equation}
\label{f!}
(\varepsilon (f_i)_! [c] )(x_i , g) = 1 \text{ whenever } g \in \langle l \rangle _H \subset G,
\quad \text{ and $0$ otherwise. }
\end{equation}

We turn to the pairing $\langle \varepsilon (f_i)_! [c] , \ch [ D ] \rangle $.
By Ponge-Wang \cite[Theorem 11.4]{PW18},
the pairing $\langle \varepsilon (f_i)_! [c] , \ch [ D ] \rangle $ is computed by taking the cap product of \eqref{f!}
and \eqref{PWCh}, and then integrate over $M ^g$.
Since every $g \in \langle l \rangle_H$ is non-trivial and $M^g$ is zero dimensional,
the integration becomes evaluation.
Hence we get
\begin{equation}
\label{DisPair}
\langle \varepsilon (f_i)_! [c] , \ch [ D ] \rangle
= \sum _{g \in \langle l \rangle _H} \Upsilon^g (R) (x_i, g)
= \sum _{g \in \langle l \rangle _H}
(-i)^{\frac{n}{2}}\det {}_\CC^{-\frac{1}{2}}(1-g^{\mathcal N}_{x_i}).
\end{equation}

Hence taking \eqref{DisPair} then summing over all $x_i$,
and taking \eqref{AB} and then summing over all $g \in \langle l \rangle _H$,
and comparing the two, we conclude that
\begin{equation}
\label{IndSum}
\sum _{i=1}^N \langle \varepsilon (f_i)_! [c] , \ch [ D ] \rangle
= \sum _{i=1}^N \sum _{g \in \langle l \rangle _H} \Upsilon^g(R)(x_i,g)
= \sum _{g \in \langle l \rangle _H} \ind _g (D). \qedhere
\end{equation}
\end{proof}

Recall that $\ind _g (D) $ is an adjoint invariant function of $g \in G$.
In particular, for any $g \in \langle l \rangle _H$ (which implies $g \in \langle l \rangle _G$),
$$ \ind _g (D) = \ind _l (D).$$
Combining with Equation \eqref{IndSum}, we recover the $l$-component
of the finite-group equivariant index from our index pairing.
\begin{cor}
We have
$$ \sum _{i=1}^N \langle \varepsilon (f_i)_! [c] , \ch [ D ] \rangle = |\langle l \rangle _H| \ind _l (D).$$
\end{cor}

In the more general case when the points with non-trivial isotropy are discrete,
but each can have a different isotropy subgroup,
a formula similar to \eqref{IndSum} can also be derived.
Now we suppose that we are given a fixed $l \in G \setminus \{ e \}$.
Our objective is to compute $\ind _l (D)$.

Again we denote by $M ^l$ the fixed point set of $l$.
We may assume $M^l \neq \emptyset$ otherwise the pairing just vanishes,
and we let $M ^l = \{ x_1 , \cdots, x_N \}$.
For each $x_i$ , we denote the isotropy subgroup of $x_i$ in $G$ by $G_i \; (:= \maG _{x_i} ^{x_i})$.
Note that $l \in G_i$ for all $i$,
hence it makes sense to define $[ c _i ] \in H P _0 (\CC G_i) \; ( = H P _0 (\maC ^\infty _c (\maG _{x_i} ^{x_i}))$ by
$$ c _i := \sum _{g \in \langle l \rangle _{G_i} } g.$$
Again we consider the pairing $\langle \varepsilon (f_i)_! [c _i] , \ch [ D ] \rangle $.

\begin{prop}
\label{EquivPairing}
We have
$$ \ind_l(D) = \sum _{i=1}^N \frac{1}{|\langle l \rangle_{G_i}|} \langle \varepsilon (f_i)_! [c _i] , \ch [ D ] \rangle .$$
\end{prop}

\begin{proof}
The computation leading to Equation \eqref{DisPair} is still valid and gives
$$ \langle \varepsilon (f_i)_! [c _i] , \ch [ D ] \rangle
= \sum _{g \in \langle l \rangle _{G_i}} \Upsilon^g (R) (x_i, g) .$$
Because the local contribution \eqref{PWIsolated} is invariant under conjugation in $G_i$,
it follows that
$$ \sum _{g \in \langle l \rangle _{G_i}} \Upsilon^g (R) (x_i, g) = |\langle l \rangle_{G_i}|\,\Upsilon^l(R)(x_i,l) .$$
On the other hand, the fixed-point formula gives
$$ \ind_l(D) = \sum_{i=1}^N \Upsilon^l (R) (x_i, l)
= \sum _{i=1}^N \frac{1}{|\langle l \rangle_{G_i}|} \langle \varepsilon (f_i)_! [c _i] , \ch [ D ] \rangle ,$$
as asserted by the proposition.
\end{proof}

%

\begin{rem}
For Morita equivalent $M \rtimes G$ and $M' \rtimes G'$ it is not obvious how the fixed point sets correspond to each other
(assuming they are discrete),
so one still cannot compare the equivariant index on these spaces by simply applying Proposition \ref{EquivPairing}.
\end{rem}

\appendix

\section{Sheaf homology of an \'etale groupoid}
In this appendix, we recall some basic facts about $\maG$-sheaves and their groupoid homology.
We shall only consider sheaves over paracompact and Hausdorff topological spaces (unlike Moerdijk).

\subsection{Basic definition of sheaves}
We begin with recalling some basic definitions of sheaves;
our exposition largely follows \cite[Chapter 10]{GallierQuaintanceBook25} and \cite[Chapter 2]{KashiwaraShapiraBook90}.

Let $M$ be a paracompact and Hausdorff topological space.
\begin{dfn}
A pre-sheaf $\maF$ is a family of local sections $\maF (U )$ for each open subset $U \subset M$.
For any open $U' \subset U$, there is a restriction map $r_{U,U'} : \maF (U ) \to \maF (U')$,
which satisfies
$r _{U, U''} = r_{U', U''} \circ r_{U, U'}$
whenever $U'' \subset U' \subset U$.

If each $\maF (U)$ is an Abelian group and all $r_{U, U'}$ are group homomorphisms, $\maF$ is called a pre-sheaf of Abelian groups.
\end{dfn}

\begin{dfn}
A pre-sheaf $\maF$ is called a sheaf if for any open $U \subset M$ and cover
$ U _\alpha $ of $U$,
\begin{enumerate}
\item
for any $s, s' \in \maF (U)$, if $r_{U, U _\alpha} (s') = r_{U, U _\alpha} (s) $
for all $\alpha$, then $s = s'$;
\item
for any collection $s _\alpha \in \maF (U _\alpha )$, such that
$$ r _{U _\alpha , U_\alpha \cap U_\beta} (s _\alpha ) = r _{U_\beta , U_\beta \cap U_\alpha} (s _\beta )
\text { for all $\alpha, \beta$},$$
there exists a unique $s \in \maF (U )$ such that
$r _{U ,U _\alpha} (s) = s _\alpha$ for all $\alpha$.
\end{enumerate}
\end{dfn}

\begin{dfn}
Let $\maF$ be a pre-sheaf and $x \in M$ be any fixed element.
Define an equivalence relation on $\maF$, where for any neighborhoods $U, U'$ of $x$,
$s \in \maF (U), s' \in \maF (U')$,
$$ s \sim s' \text{ if } r _{U ,U''} (s) = r _{U', U''}(s') \text{ on some } U'' \subset U \cap U'.$$
The set of equivalence classes is called the stalk of $\maF$ at $x$, denoted $\maF _x$.
In other words
$$\maF _x := \lim _\to \{ \maF (U) : x \in U \}.$$
Each germ is an Abelian group if $\maF$ is a pre-sheaf of Abelian groups.
We also denote $s _x := s / \sim$ and call it the germ of $s$ at $x$.
\end{dfn}

\begin{dfn}
Let $\maF$ be a pre-sheaf over $M$. The \'etale space of $\maF$ is the disjoint union
$$\maE _\maF := \coprod _{x \in M} \maF _x,$$
endowed with the topology generated by all subsets of the form
$$\{ s_x : s \in \maF (U), x \in U \}.$$
There is a natural projection from $\maE _\maF$ to $M$.
\end{dfn}

\begin{dfn}
Let $\maF$ be a pre-sheaf. The sheafification of $\maF$ is the sheaf $\maF ^+$
defined by
\begin{align*}
\maF ^+ & (U) := \Big\{ s : U \to \maE _\maF \text{ such that for all $x \in U, s(x) \in \maF _x$, and there exists a} \\
& \text{neighborhood $U' \subset U$ and $s' \in \maF (U')$ satisfying $s(x') = s' _{x'}$ for all $x' \in U'$}. \Big\}
\end{align*}
If $\maF$ is already a sheaf then $\maF ^+$ is canonically isomorphic to $\maF$.
Hence $\maF (U)$ is just the group of all continuous sections of $\maE _\maF$ over $U$.
\end{dfn}

\begin{dfn}
For any $s \in \maF (M)$, the support of $s$ is the complement in $M$ of the subset
$$ \{ x \in M :
\text{ there exists a neighborhood $U$ of $x$ such that } r _{M,U} (s) = 0 \}.$$
We say that $s$ is compactly supported if the support of $s$ is a compact set,
and denote the group of all such sections by $\Gamma _c (M, \maF)$.
\end{dfn}

\begin{dfn}
\cite[Definition 2.3.1 (ii)]{Moerdijk;EtaleHomology}
Let $M, M'$ be topological spaces, $\varphi : M' \to M$ be a continuous map.
Given a sheaf $\maF$ over $M$, the inverse image $\varphi^* \maF$ is the sheafification of the pre-sheaf over $M'$ by
$$ \lim _\rightarrow \{ \maF (U' ) : U' \supset \varphi (U)\}$$ open.
Note that the stalk of $\varphi ^* \maF$ at $y \in M'$ is $\maF _{\varphi (y)}$.
\end{dfn}

\begin{dfn}
\cite[Definition 2.3.1 (i), Equation (2.5.1)]{Moerdijk;EtaleHomology} Let $M, M'$ be topological spaces,
$\varphi : M' \to M$ be a continuous map. Given a sheaf $\maF'$ over $\maM'$, the direct
image sheaf is the pre-sheaf
$$\varphi _* \maF' (U ) := \maF' (\varphi^{-1} U)$$
(which is a sheaf).

Suppose furthermore that $M ,M'$ are both locally compact, define the compactly supported direct image sheaf to be
$$ \varphi _! \maF' (U) := \{ s \in \maF' (\varphi^{-1} U ) : \varphi \vert_{\supp(s)} \text{ is proper} \}.$$
\end{dfn}
Observe that one has canonical isomorphism
$$ \Gamma _c (M' , \maF') \cong \Gamma _c (\{p \}, p_! \maF') \cong \Gamma _c (M, \varphi _! \maF'),$$
where $p$ is the map into the one-point space. Moreover, the stalk of $\varphi_! \maF'$ is
isomorphic to
$$\Gamma _c (\varphi^{-1} (x), \maF' \vert _{\varphi^{-1}(x)}).$$

\begin{ex}
\label{BackForwardSheaf}
Let $M, M'$ be locally compact, Hausdorff topological spaces, and $\varphi : M' \to M$ be a local homeomorphism.
Then for any $x \in M$, $\varphi ^{-1} (x)$ is a discrete subset of $M'$.
Let $\maF$ be any sheaf of Abelian groups on $M$.
From definition, the stalk of $\varphi _! \varphi ^* \maF$ at $x \in M$ is
$$ (\varphi _! \varphi ^* \maF) _x = \Gamma_c (\varphi^{-1} (x), (\varphi^* \maF) \vert _{\varphi^{-1}(x)})
= \bigoplus _{y \in \varphi ^{-1} (x)} (\varphi^* \maF)_y = \bigoplus _{y \in \varphi ^{-1} (x)} \maF _x .$$
In particular there is a canonical map from $\Gamma _c (M, \varphi _! \varphi ^* \maF) $ to $\Gamma _c (M, \maF)$
defined by summation over the fibers.
\end{ex}

\begin{dfn}
We say that $\maF$ is c-soft if for any compact subset $K \subset M$,
the restriction morphism from $\maF (M)$ to $\maF (K)$ is surjective,
where $\maF(K):=j^* \maF(M)$ and $j:K \hookrightarrow M$ is the inclusion map.
\end{dfn}

In practice, we are mostly interested in the following examples.
\begin{ex}
Let $M$ be a manifold. For each open subset $U \subset M$, let $\maF (U)$ be the
vector space of smooth functions on $U$.
This defines a sheaf, known as the sheaf of germs of smooth functions on $M$, $C^\infty _M$.
Similarly, we may define the sheaf of germs of differential forms $\Omega^\bullet _M$.
These sheaves are fine and therefore c-soft.
\end{ex}

\subsection{$\maG$-sheaves and homology}
In this section we only consider Lie groupoids that are \'etale.

\begin{dfn}
A $\maG$-sheaf is a sheaf $\maF$ over $\maG^{(0)}$,
such that the \'etale space $\maE _\maF$ is
endowed with a continuous right action of $\maG$.
In other words, for each element $\gamma \in \maG$
one is given a group isomorphism $\gamma : \maF _{\bt(\gamma)} \to \maF _{\bs(\gamma)}$, that is associative,
is identity if $\gamma$ is a unit, and induces an isomorphism $\maF (\bt(U)) \to \maF (\bs(U))$ for any bisection $U$.
\end{dfn}

\begin{ex}
\label{SmoothGerms}
Let $\maG$ be an \'etale groupoid, and consider the sheaf $C^\infty _{\maG ^{(0)}}$ over $\maG ^{(0)}$.
For any $\gamma \in \maG ^{(1)}$ there exists a bisection $S$ containing $\gamma$.
Then $\bt \circ (\bs |_S)^{-1}$ maps a neighborhood of $\bs (\gamma)$ to $\bt (\gamma)$,
hence pull back by $\bt \circ (\bs |_S)^{-1}$ defines a map from the germ of any smooth function at $\bt (\gamma)$ to $\bs (\gamma)$,
that is independent of $S$.
\end{ex}

Let $\maG$ be a groupoid. We denote
$$\maG ^{(n)} := \{(\gamma _1 , \cdots, \gamma_n ) \in \maG^{(1)} \times \cdots \times \maG^{(1)} :
\bs(\gamma_i) = \bt(\gamma _{i+1}), i = 1, \cdots, n-1 \}.$$

Recall that $\maG ^{(\bullet)}$ is endowed with a simplicial manifold structure,
with the face maps $\partial _i : \maG^{(n)} \to \maG ^{(n-1)}$
\begin{align}
\nonumber
\partial _0 (\gamma _1 , \cdots , \gamma _n ) :=& (\gamma _2, \cdots, \gamma _n) \\
\partial _i (\gamma _1 , \cdots , \gamma _n ) :=& (\gamma _1 , \cdots \gamma _i \gamma _{i+1} \cdots , \gamma _n ),
\quad 1 \leq i \leq n-1 \\ \nonumber
\partial _n (\gamma _1 , \cdots , \gamma _n ) :=& (\gamma _1, \cdots, \gamma _{n-1}).
\end{align}

Let $\maA$ be any c-soft $\maG$-sheaf.
We define $\bt_n : \maG ^{(n)} \to \maG^{(0)}$
\begin{equation}
\label{tn}
\bt_n (\gamma _1 , \cdots , \gamma _n ) := \bt ( \gamma _1 ),
\end{equation}
and let
$$\maA _n := \bt _n ^* \maA .$$
\begin{dfn}
\label{BarComplex}
We denote
$$ B_n (\maG , \maA) := \Gamma _c (\maG ^{(n)}, \maA _n). $$
Define $\partial_i : B_n (\maG, \maA) \to B_{n-1} (\maG, \maA)$
by summation along the fiber of $\partial_i : \maG^{(n)} \to \maG ^{(n-1)}$ (see Example \ref{BackForwardSheaf}).
Explicitly, for any $c \in B_n (\maG, \maA)$,
\begin{align*}
\label {FaceMaps}
\partial _0 c (\gamma _1 , \cdots , \gamma _{n-1} )
:=& \sum _{\bs (\gamma ) = \bt (\gamma _1)} c (\gamma , \gamma _1 , \cdots , \gamma _{n-1} ) \cdot \gamma \\
\partial _i c (\gamma _1 , \cdots , \gamma _{n-1} )
:=& \sum _{\bt (\gamma ) = \bt (\gamma _{i})} c (\gamma _1 , \cdots , \gamma, \gamma ^{-1} \gamma _i, \cdots, \gamma _{n-1}) \\
\partial _n c (\gamma _1 , \cdots , \gamma _{n-1} )
:=& \sum _{\bt (\gamma ) = \bs (\gamma _{n-1})} c (\gamma _1 , \cdots , \gamma _{n-1}, \gamma).
\end{align*}
Then $\partial := \sum _{i=0}^n (-1)^i \partial _i$ turns $B _n (\maG , \maA)$ into a chain complex,
called the bar complex, its homology is known as the (compactly supported) groupoid homology with coefficient $\maA$
$$ H _\bullet (\maG, \maA ).$$
\end{dfn}

To define the groupoid homology of general sheaves, we consider c-soft resolutions.

\begin{dfn}
\label{ChainCyclic}
Let
$$ \cdots \to \maF ^\bullet \xrightarrow{d} \maF ^{\bullet - 1} \to \cdots $$
be a bounded below complex of c-soft $\maG$-sheaves, the groupoid homology
$$ H _\bullet (\maG, \maF ^\bullet )$$
is defined to be the homology of the total complex
$$ H _\bullet \Big(\bigoplus _{p+q=\bullet} B_p (\maG , \maF ^q), \partial + d \Big) .$$
\end{dfn}

\begin{dfn}
Let $\maA$ be a $\maG$-sheaf and suppose we have a resolution of $\maA$ by c-soft sheaves,
i.e., a bounded below complex of $\maG$-sheaves
$$ 0 \to \maA \to \maF ^0 \xrightarrow{d} \maF ^{-1} \to \cdots, $$
where $\maF ^{-k}$ are c-soft, and the complex is exact at each stalk.
Then we define the groupoid homology with coefficient $\maA$ by
$$ H _\bullet (\maG, \maA) := H _\bullet \Big(\bigoplus _{p+q=\bullet} B_p (\maG , \maF ^q), \partial + d \Big) .$$
\end{dfn}
%

Here, we recall the functoriality of groupoid homology \cite[Section 3.7]{Moerdijk;EtaleHomology}.
\begin{lem}
\label{HomFunct}
Suppose that $\varphi : \maK \to \maG$ is proper, then for any $\maG$-sheaf $\maA$ one obtains homomorphisms
$$ \Gamma_c (\maG ^{(n)} , \maA _n) \to \Gamma_c (\maK ^{(n)} , \varphi ^* \maA _n) $$
by pullback, and hence a homomorphism:
$$ \varphi ^* : H _n (\maG , \maA) \to H_n (\maK , \varphi^* \maA) .$$

Suppose $\varphi : \maK \to \maG$ is \'etale (i.e., each $\varphi _n : \maK ^{(n)} \to \maG ^{(n)}$ is a local diffeomorphism).
Let $\maS$ be a c-soft $\maG$-sheaf.
Summation along the fibers defines a homomorphism
$$(\varphi _n )_! \bt _n^* (\varphi _n)^* \maS = (\varphi _n)_! \varphi _n^* \bt _n ^* \maS \to \bt _n ^* \maS $$
and hence homomorphisms
$$ \Gamma_c (\maK ^{(n)} , \varphi ^* \maS _n) \to \Gamma_c (\maG ^{(n)} , \maS _n) $$
that commute with the face operators. That in turn induces homomorphisms
$$ \varphi _* : H_n (\maK , \varphi^* \maA) \to H_n (\maG , \maA) $$
for any $\maG$-sheaf $\maA$ by taking c-soft resolution.
\end{lem}

\begin{ex}
\label{ThomEx}
Let $G \subset G L (k, \RR)$ be a discrete subgroup acting on $\RR ^k$ linearly, and preserving the orientation.
Let $ \maG := \RR ^k \rtimes G$ be the action groupoid.
Hence elements of $\maG ^{(n)}$ are
$$ ((x, g_1), (x g_1 , g_2), \cdots (x g_1 \cdots g_{n-1}, g_n)),$$
which we identify with
$$ (x, g_1, \cdots , g_n) \in \RR ^k \times G^n .$$

To compute $H _\bullet (\maG , \RR)$, we consider the double complex
$$
\begin{CD}
@. \vdots @. \vdots \\
@. @VVV @VVV \\
\cdots @<<< \Gamma_c (\maG ^{(1)} , \bt_1^* \Omega ^1) @<{d _{\DR}}<< \Gamma_c (\maG ^{(1)} , \bt_1^* \Omega ^0) \\
@. @V{\partial}VV @V{\partial}VV \\
\cdots @<<< \Gamma_c (\maG ^{(0)} , \Omega ^1) @<{d _{\DR}}<< \Gamma_c (\maG ^{(0)} , \Omega ^0)
\end{CD}.
$$

Consider the rows.
An element in $\Gamma_c (\maG ^{(p)} , \bt_p^* \Omega ^q)$
is just a family of compactly supported $q$-forms on the connected components
$\RR ^k \times {(g_1, \cdots , g_n)}$,
and $d_{\DR}$ is just the de Rham differential on each component.
It is well known, by the Poincar\'e lemma with compact support, that
$$ H^k (\Omega _c ^\bullet (\RR ^k) , d _{\DR}) \cong \CC, \quad H ^i (\Omega _c ^\bullet (\RR ^k), d _{\DR}) \cong 0
\text{ for $0 \leq i < k$}. $$
Moreover, integrating the top form gives a canonical isomorphism that is invariant under orientation-preserving diffeomorphisms.
It follows that integration on each connected component
$\RR ^k \times {(g_1, \cdots , g_n)}$
defines a map
$$\maI : \Gamma_c (\maG ^{(p)} , \bt_p^* \Omega ^k) \to \Gamma _c (G^p, \CC),$$
which commutes with $\partial$ and is a quasi-isomorphism of the rows.
It follows that the total complex collapses to the group homology, i.e. we arrive at
$$ H _{\bullet - k} (\maG , \CC) \cong H _\bullet (G, \CC) .$$
The inverse $\maI ^{-1} (u)$ is represented by an element in $\Gamma_c (\maG^{(p)} , \bt_p^* \Omega ^q)$
whose top de Rham degree part integrates to $u (\gamma _1, \cdots, \gamma _p)$.

Note that one can replace $\RR ^k$ by any $G$-invariant contractible open subset and the arguments above hold without change.
\end{ex}

\begin{rem}
The example above can be generalized.
First, let $U \subset \RR ^k$ be any star shaped neighborhood of $0 \in \RR ^k$.
Then
$$ \maG _U := \bs ^{-1} U \cap \bt ^{-1} U \subset \maG $$
is a subgroupoid.
Since for any $g_1, \cdots , g_p \in G$,
$$ g_1 \cdot U \cap \cdots \cap g_i \cdot U $$
is also star shaped,
all the arguments above works without change, and we have
$$ H _{\bullet - k} (\maG _U , \CC) \cong H _\bullet (G, \CC) .$$
It seems that if the groupoid in Section 4.1 is not locally a action groupoid (condition (3)),
one can consider $\maG _U $ instead.
\end{rem}



\end{document}